\documentclass[11pt,twoside]{article}
\usepackage{latexsym}
\input{amssym.def}
\input{amssym}
\newfont{\fra}{eufm10 scaled 1095}
\newfont{\Bb}{msbm10 scaled 1095}
\newfont{\Bbg}{msbm10 scaled 1280}
\newcommand\CC{{\mbox{\Bb C}}}
\newcommand\RR{{\mbox{\Bb R}}}
\newcommand\NN{{\mbox{\Bb N}}}
\newcommand\ZZ{{\mbox{\Bb Z}}}

\newcommand\fg{{\frak{g}}}
\newcommand\fh{{\frak h}}
\newcommand\fri{{\frak i}}
\newcommand\fj{{\frak j}}
\newcommand\fl{{\frak l}}
\newcommand\fm{{\frak m}}
\newcommand\fn{{\frak n}}
\newcommand\fa{{\frak a}}
\newcommand\fb{{\frak b}}
\newcommand\fd{{\frak d}}

\newcommand\fk{{\frak k}}
\newcommand\fr{{\frak r}}

\newcommand\fz{{\frak z}}

\newcommand\cZ{{\cal Z}}
\newcommand\cH{{\cal H}}

\newcommand\cA{{\cal A}}

\newcommand\ph{\varphi}

\newcommand\tha{{\theta_{\fa}}}
\newcommand\thl{{\theta_{\fl}}}

\newcommand{\fsl}{\mathop{{\frak s \frak l}}}
\newcommand{\fsu}{\mathop{{\frak s \frak u}}}

\newcommand{\so}{\mathop{{\frak s \frak o}}}

\newcommand{\Aut}{\mathop{{\rm Aut}}}
\newcommand{\GL}{\mathop{{\it GL}}}

\newcommand{\Hom}{\mathop{{\rm Hom}}}
\newcommand{\Id}{\mathop{{\rm Id}}}

\newcommand{\ad}{{\rm ad}}

\newcommand{\Ad}{\mathop{{\rm Ad}}}

\newcommand{\Ker}{\mathop{{\rm ker}}}

\newcommand{\Span}{\mathop{{\rm span}}}
\newcommand{\mod}{\mathop{{\rm mod}}}
\newcommand{\proj}{{\rm pr}}
\newcommand\ip{\mbox{$\langle\cdot \,,\cdot \rangle$}}
\newcommand\ipa{{\langle\cdot \,,\cdot \rangle_\fa}}
\newcommand\lb{{[\cdot\,,\cdot]}}
\newcommand\dd{\fd_{\alpha,\gamma}(\fl,\thl,\fa)}
\newcommand\proof{{\sl Proof. }}
\newcommand{\qed}{\hspace*{\fill}\hbox{$\Box$}\vspace{2ex}}
\newtheorem{theo}{Theorem}[section]
\newtheorem{pr}{Proposition}[section]
\newtheorem{de}{Definition}[section]

\newtheorem{re}{Remark}[section]
\newtheorem{co}{Corollary}[section]
\newtheorem{lm}{Lemma}[section]
\begin{document}
\title{On the structure of pseudo-Riemannian symmetric spaces}
\author{Ines Kath and Martin Olbrich}
\maketitle
\begin{abstract}
\noindent
Following our approach to metric Lie algebras developed in a previous
paper we propose a way of understanding pseudo-Riemannian symmetric spaces
which are not semi-simple. We introduce cohomology sets (called
quadratic cohomology) associated with orthogonal
modules of Lie algebras with involution. Then we construct a functorial assignment
which sends a pseudo-Riemannian symmetric space $M$ to a triple consisting of\\
(i) a Lie algebra with involution (of dimension much smaller than the dimension\\
\makebox{$\quad   $} of the transvection group of $M$),\\
(ii) a semi-simple orthogonal module of the Lie algebra with involution, and\\
(iii) a quadratic cohomology class of this module.\\
That leads to a classification scheme of indecomposable non-simple pseudo-Rieman\-nian symmetric spaces. In addition, we obtain a full
classification of symmetric spaces of index 2 (thereby completing and correcting in part earlier classification results due to Cahen/Parker and
Neukirchner).
\end{abstract}
\tableofcontents

\section{Introduction}

In the present paper we will develop a classification scheme for pseudo-Riemannian symmetric spaces. A (pseudo-)Riemannian symmetric space is a connected (pseudo-)\linebreak Riemannian manifold $(M,g)$ satisfying the following symmetry condition. For each point $x\in  M$ there is an involutive isometry $\theta_x$ of $(M,g)$ such that $x$ is an isolated fixed point of $\theta_x$. Clearly, in this case the ``symmetry'' $\theta_x$ is uniquely defined. In a normal neighbourhood of $x$ it coincides with the geodesic reflection at $x$. Riemannian symmetric spaces are well understood. They were classified by \'E.~Cartan. Contrary to that the classification in the pseudo-Riemannian case is much more involved.

One can reformulate the classification problem for pseudo-Riemannian symmetric spaces in terms of purely algebraic objects. In order to do so we assign to each pseudo-Riemannian symmetric space $M$ a distinguished subgroup $G$ of its group of isometries, the transvection group of $M$ (sometimes also called
group of displacements) which still acts transitively on $M$. For a fixed base point $x_0\in M$ the involution $\theta_{x_0}$ induces an involution $\theta$ on the Lie algebra $\fg$ of $G$ and the scalar product $g_{x_0}$ on the tangent space $T_{x_0}M$ at $x_0$ induces a (non-degenerate) scalar product $\ip$ on $\fg$, which is $\fg$- and $\theta$-invariant. Moreover, the decomposition $\fg=\fg_+\oplus \fg_-$ of $\fg$ into the $\pm 1$-eigenspaces of $\theta$ satisfies $[\fg_-,\fg_-]=\fg_+$ (see Section~\ref{sym} for a more exact explanation of these facts). A triple $(\fg,\theta,\ip)$ which consists of a Lie algebra $\fg$, a scalar product $\ip$ and an involution $\theta$ satisfying these conditions will be called a symmetric triple. This gives us a one-to-one correspondence between simply-connected pseudo-Riemannian symmetric spaces and symmetric triples. Hence the classification of simply-connected pseudo-Riemannian symmetric spaces up to isometry is equivalent to the
classification of symmetric triples up to isomorphism.

A classification of all semi-simple symmetric triples, i.e.~of symmetric triples $(\fg,\ip,\theta)$ with semi-simple
$\fg$, has been already obtained by M.~Berger \cite{B}. In particular, this includes all symmetric triples which correspond to non-flat irreducible pseudo-Riemannian symmetric spaces. However, contrary to the Riemannian case a pseudo-Riemannian symmetric space in general does not decompose into a product of irreducible symmetric spaces. Thus Berger's result does not solve our classification problem. In fact, symmetric spaces with semi-simple transvection group constitute only a very small part of all pseudo-Riemannian symmetric spaces. The next step was done by M.~Cahen and N.~Wallach \cite{CW} who classified all symmetric triples which correspond to Lorentzian symmetric spaces, i.e.~to pseudo-Riemannian symmetric spaces with a metric of index~1. For the case of index~2
M.\,Cahen and M.\,Parker gave classification results in \cite{CP1} and \cite{CP}. However the results in \cite{CP1} are not complete and also the revised version of these results by Th.~Neukircher in \cite{N} is not quite correct (see also the introduction to Section~\ref{S7}). For higher index the classification problem is unsolved. Note however, that there are recent classification results
of pseudo-Riemannian symmetric spaces admitting very special additional
geometric structures ((para-) hyper-K\"ahler \cite{AC1},\cite{ABCV}, (para-) quaternionic K\"ahler \cite{AC2}).

Looking at the results for index 1 and 2 it becomes rather obvious that one cannot expect to get a list of all symmetric triples for arbitrary index. Therefore the aim will be to find a reasonable description of the structure of symmetric triples which will lead to a nice description of the moduli space of isomorphism classes of all symmetric triples.
For obvious reasons one is mainly interested in indecomposable symmetric triples, i.e. in those which are not the non-trivial direct sum of symmetric triples. Moreover it suffices to consider only those symmetric triples $(\fg,\theta,\ip)$ for which $\fg$ does not contain semi-simple ideals. Indeed, if $\fg$ has a semi-simple ideal, then it also has a $\theta$-invariant semi-simple ideal $\fh$. The restriction of $\ip$ to any semi-simple ideal is non-degenerate and therefore the symmetric triple decomposes into the direct sum $(\fh,\theta|_{\fh},\ip_{\fh})\oplus(\fh^\perp,\theta|_{\fh^\perp},\ip_{\fh^\perp})$
of two symmetric triples the first one being in Berger's classification list mentioned above.

The structure theory for symmetric spaces which we will develop in this paper is completely parallel to our theory of metric Lie algebras presented in \cite{KO2}.
Let us describe the main ideas and results in the present context of symmetric triples. For each symmetric triple $(\fg,\theta,\ip)$ without semi-simple ideals there is a canonically defined $\theta$-invariant isotropic
ideal $\fri(\fg)\subset \fg$ such that $\fri(\fg)^\perp/\fri(\fg)$ is abelian
(see Section \ref{poseidon}, Formula~(\ref{herakles})).
Furthermore, for any $\theta$-invariant isotropic
ideal $\fri\subset \fg$ such that $\fri^\perp/\fri$ is abelian the following holds. Put
$\fl=\fg/\fri^\perp$ and $\fa=\fri^\perp/\fri$. The involution $\theta$ induces involutions $\theta_\fl$ on $\fl$ and $\tha$ on $\fa$. Then $\fa$ inherits
an inner product from $\fg$ and an $\fl$-action being compatible with this inner product and the involutions $\tha$ and $\thl$, i.e., it inherits the structure of an orthogonal $(\fl,\thl)$-module. Moreover, $\fri\cong \fl^*$ as an $(\fl,\thl)$-module, and $\fg$ can be represented
as the result of two subsequent extensions of Lie algebras with involution
\begin{equation}
0\rightarrow \fa{\longrightarrow}\tilde\fg
{\longrightarrow} \fl\rightarrow 0\ ,\qquad
0\rightarrow \fl^*{\longrightarrow}\fg
{\longrightarrow} \tilde\fg\rightarrow 0\ ,\label{lena}
\end{equation}
where $\tilde\fg=\fg/\fri$.
Vice versa, given a Lie algebra with involution $(\fl,\thl)$, an orthogonal $(\fl,\thl)$-module
$\fa$ and two extensions as in (\ref{lena}) which in addition satisfy
certain natural compatibility conditions, the resulting Lie algebra $\fg$
has a distinguished invariant inner product and an isometric involution. This construction will be formalised into the notion of a quadratic extension of $(\fl,\thl)$ by
an orthogonal $(\fl,\thl)$-module $\fa$ in Subsection \ref{S41}. In particular, there is a natural
equivalence relation on the set of quadratic extensions of $(\fl,\thl)$ by $\fa$.
We will call such a quadratic extension of a Lie algebra with involution $(\fl,\thl)$ admissible if the resulting metric Lie algebra with involution
$(\fg,\thl,\ip)$ is a symmetric triple and the image of $\fl^*$ in $\fg$ coincides with the canonical ideal $\fri(\fg)\subset \fg$ (the crucial condition).
Admissibility implies in particular that the $(\fl,\thl)$-module $\fa$ is semi-simple.

Summarising up to here we have established that any symmetric triple
without semi-simple ideals has the structure of an admissible quadratic extension of a Lie algebra
with involution $(\fl,\thl)$ by an orthogonal $(\fl,\thl)$-module $\fa$
in a canonical way.
Geometrically, the ideal $\fri(\fg)^\perp\subset\fg$ defines a
foliation of the pseudo-Riemannian space~$M$ invariant under all transvections. Its leaves are coisotropic symmetric subspaces of $M$ which are flat. This foliation is actually a fibration over an affine symmetric space~$N$ (without metric). The Lie algebra $\fl$ is then
a central extension of the Lie algebra of the transvection group of $N$.
The theory discussed so far is essentially a formalisation of ideas of
L.~Berard-Bergery \cite{BB1},\cite{BB2} (unpublished).

This formalisation
allows us to proceed further. First, equivalence classes of quadratic extensions
of $(\fl,\thl)$ by an orthogonal module $\fa$ are conveniently described by a certain cohomology
set $\cH_Q^2(\fl,\thl,\fa)$ (which is the quotient of a real algebraic variety
by an algebraic group action).
We introduce these quadratic cohomology sets for orthogonal modules of
Lie algebras with involution in Section \ref{grishkov}. The relation to quadratic
extensions is established in Section \ref{leda}. In particular,
Subsection~\ref{neptun} contains our standard construction which produces a metric Lie algebra with involution $(\fg,\theta,\ip)$ out of a Lie algebra with involution $(\fl,\thl)$, an orthogonal $(\fl,\thl)$-module $\fa$, and a cocycle representing
an element in $\cH_Q^2(\fl,\thl,\fa)$.

Secondly, for semi-simple modules $\fa$, we describe the subset $\cH_Q^2(\fl,\thl,\fa)_\sharp\subset\cH_Q^2(\fl,\thl,\fa)$ corresponding to
admissible quadratic extensions (Section \ref{poseidon}).
We call a Lie algebra with involution $(\fl,\thl)$ admissible if it possesses
a semi-simple orthogonal module $\fa$ such that
$\cH_Q^2(\fl,\thl,\fa)_\sharp\ne \emptyset$. As examples show this condition is  rather strong. In particular, not every affine symmetric space can appear
as the base $N$ of the special fibration described
above associated with some pseudo-Riemannian symmetric space $M$.

And thirdly, in Section \ref{quark} we establish that the correspondence
$$ \{ \mbox{symmetric triples } (\fg,\theta,\ip)\}\Longrightarrow \left\{\mbox{quadruples } \left(\fl,\thl,\fa,\varphi\in
\cH_Q^2(\fl,\thl,\fa)_\sharp\right)\right\}\ , $$
which sends each symmetric triple to the datum which defines the equivalence class of its associated canonical quadratic extension descends to a bijection
of isomorphism classes. Here two quadruples $(\fl_i,\theta_{\fl_i},\fa_i,\varphi_i)$,
$i=1,2$, are called to be isomorphic if there is an isomorphism of triples
(for the formal
definition of this notion see the end of Section~\ref{grishkov})
$$T:(\fl_1,\theta_{\fl_1},\fa_1)\rightarrow (\fl_2,\theta_{\fl_2},\fa_2)\  \mbox{ such that }\  T^*\varphi_2=\varphi_1\ .$$
Let us express this result in a slightly different
way. We consider the automorphism group $G_{\fl,\theta_\fl,\fa}$ of a fixed triple
$(\fl,\thl,\fa)$. It is a subgroup of $\Aut(\fl,\thl)\times O(\fa)^\tha$.
Then the moduli space of all symmetric triples without semi-simple ideals
can be identified with the union of quotient spaces
\begin{equation}\label{class}
 \coprod_{(\fl,\theta_\fl,\fa)} \cH^2_Q(\fl,\thl,\fa)_\sharp/G_{\fl,\theta_\fl,\fa}\ ,
\end{equation}
where the union is taken over a set of representatives of isomorphism
classes of triples $(\fl,\thl,\fa)$ consisting of an admissible Lie algebra
with involution $(\fl,\thl)$ and a semi-simple orthogonal $(\fl,\thl)$-module $\fa$.
As already mentioned, we are mainly interested in indecomposable symmetric
triples. Definition \ref{kosel} together with Proposition \ref{hermlin} provides
a manageable description of the subset $\cH^2_Q(\fl,\thl,\fa)_0\subset \cH^2_Q(\fl,\thl,\fa)_\sharp$ corresponding
to indecomposable symmetric triples. Therefore the moduli space
of non-semi-simple indecomposable symmetric triples is given by the following subspace of (\ref{class})
\begin{equation}\label{class0}
 \coprod_{(\fl,\theta_\fl,\fa)} \cH^2_Q(\fl,\thl,\fa)_0/G_{\fl,\theta_\fl,\fa}\ .
\end{equation}
This is the first main theorem of the present paper
(see Theorem \ref{Pwieder}).

In order to approach a true classification one has to evaluate (\ref{class0})
further. By construction, the dimension of the symmetric space $N$ corresponding to $(\fl,\thl)$ is bounded by the index of the symmetric
triple $(\fg,\theta,\ip)$, i.e.  the index of $\ip |_{\fg_-}$.
For this reason the bijection (\ref{class0}) is extremely useful for a concrete classification of symmetric triples with
small index, where all the ingredients of (\ref{class0}) can be explicitely computed. The situation concerning the classification results in index 2
which we discussed at the beginning of this introduction motivated us to redo
this classification using (\ref{class0}). This is the content of
Section \ref{S7}. Theorem \ref{Klasse} which is the second main theorem of the
present paper
provides a full classification of indecomposable
symmetric triples
of index $2$ which are not semi-simple. Moreover, we determine the
surprisingly small subset of them which belongs to pseudo-Hermitian symmetric spaces (Corollary \ref{kahl}).

Note that the holonomy representation of a symmetric space with
associated symmetric triple $(\fg,\theta,\ip)$ is given by the
adjoint representation of $\fg_+$ on $\fg_-$. Therefore the structure
results of the present paper imply structure results
for (indecomposable, non-irreducible) holonomy representations of pseudo-Riemannian symmetric
spaces. E.g., each such holonomy representation has a distinguished
invariant isotropic subspace given by $\fri(\fg)\cap\fg_-$.
As already observed by L.~Berard-Bergery \cite{BB1} a similar
canonical isotropic subspace can be assigned to any pseudo-orthogonal representation.
It should be possible to exploit this fact
in order to uncover new structural results
for holonomy representations of non-symmetric
pseudo-Riemannian manifolds as well.

\section{Symmetric pairs and triples}\label{sym}

In the following short survey of basic facts on symmetric spaces we will fix
some notation used in this paper and clarify the relations between the defined
objects. For a detailed introduction to symmetric spaces see e.\,g.\ \cite{CP},
\cite{H}, \cite{KN}.

Let $(M,\nabla)$ be a manifold with an affine connection. $(M,\nabla)$ is
called an {\it affine symmetric space} if for each $x\in M$ there is an
involutive affine transformation $\theta_x$ of
$(M,\nabla)$ such that $x$ is an isolated fixed point of $\theta_x$.

For an affine symmetric space $(M,\nabla)$ the transvection group $G$
defined by
$$G:=\langle \theta_x\circ \theta_y\mid x,y\in M\rangle.$$
acts transitively on $M$. Let us fix a point $x_0\in M$. The
involutive affine transformation $\theta_{x_0}$
acts by conjugation on $G$. We will denote this conjugation by
$\theta$. It induces an involutive automorphism of the Lie
algebra $\fg$ of $G$, which we also denote by $\theta$. We will call
$(\fg,\theta)$ {\it symmetric pair associated with} $(M,\nabla)$. Let $\fg_+$
and $\fg_-$ be the eigenspaces of $\theta$ with eigenvalues $1$ and
$-1$, respectively. Then we have
\begin{equation}\label{E0}
[\fg_-,\fg_-]=\fg_+,\qquad [\fg_+,\fg_-]\subset\fg_- .
\end{equation}
Moreover, $\fg_+$ acts faithfully on $\fg_-$.
Let $G_+\subset G$ be the stabilizer of $x_0\in M$ with respect to the
action of $G$ on $M$. Then $G_0^\theta\subset G_+\subset G^\theta$
holds. In particular, $\fg_+$ is the Lie algebra of $G_+$ whereas
$\fg_-$ can be identified with the tangent space $T_{x_0}M$ of $M$ in
$x_0$. By (\ref{E0}) the homogeneous space $G/G_+$ is reductive. In
particular, the decomposition $\fg=\fg_{+}\oplus\fg_{-}$ defines
a canonical invariant connection $\nabla ^0$. There is an affine
diffeomorphism from $(M,\nabla)$ to $(G/G_+,\nabla ^0)$.

A  {\it Lie algebra with involution} is a pair $(\fl,\theta_\fl)$ consisting
of a finite-dimensional Lie algebra $\fl$ and an involutive automorphism $\theta_\fl$ of 
$\fl$.

Let $(\fl,\thl)$ be a  Lie algebra with involution
and induced eigenspace decomposition $\fl=\fl_+\oplus\fl_-$ such that
\begin{itemize}
\item[{\rm (S1)}] $[\fl_-,\fl_-]=\fl_+$,
\item[{\rm (S2)}] $\fl_+$ acts faithfully on $\fl_-$,
\end{itemize}
then there is a uniquely determined (up to
isomorphism) simply connected affine symmetric space $(M,\nabla)$,
such that $(\fl,\theta_\fl)$ is the symmetric pair associated with
$(M,\nabla)$. Therefore we will call any Lie algebra with involution
$(\fl,\thl)$ satisfying Conditions  (S1) and (S2) a {\it symmetric pair}.

Now let $(M,g)$ be a pseudo-Riemannian manifold and denote by
$\nabla ^{LC}$ the Levi-Civita
connection on $M$. $(M,g)$ is called a {\it pseudo-Riemannian symmetric space}
if $(M,\nabla ^{LC})$ is an affine symmetric space.
In this case the affine transformations $\theta_x$ are isometries and
the transvection group $G$ is a subgroup of the isometry group.

The metric $g$ on $M$ induces a $\fg_+$-invariant (in general
indefinite) scalar product $\ip_-$ on $\fg_-$. We will call the triple
$(\fg,\theta,\ip_-)$ {\it symmetric triple associated with} $(M,g)$. It has
the following properties:
\begin{itemize}
\item[{\rm (S1)}]$[\fg_-,\fg_-]=\fg_+$,
\item[{\rm (S2)}]$\fg_+$ acts faithfully on $\fg_-$,
\item[{\rm (S3)}] $\ip_-$ has an extension to a $\theta$-invariant and
  $\fg$-invariant (non-degenerate) scalar product $\ip$ on $\fg$.
\end{itemize}
Vice versa, let $\fg$ be a finite-dimensional Lie algebra with
involutive automorphism $\theta$ and let $\fg=\fg_+\oplus\fg_-$ be the
decomposition into eigenspaces. Let $\ip_-$ be a $\fg_+$-invariant
scalar product on $\fg_-$ such that $(\fg,\theta,\ip_-)$ satisfies two of the
three Conditions
(S1), (S2) and (S3), then also the third condition
is satisfied and
there exists a uniquely defined simply connected symmetric space
$(M,g)$ such that $(\fg,\theta,\ip_-)$ is the symmetric triple
associated with $(M,g)$.
Therefore we will call every triple $(\fg,\theta,\ip_-)$ consisting of
a finite-dimensional Lie  algebra $\fg$, an involutive automorphism $\theta$ of $\fg$ 
and
an $\fg_+$-invariant
scalar product $\ip_-$ on $\fg_-$ which satisfies (two of) the Conditions
(S1), (S2) and (S3) a {\it symmetric triple}.
For a symmetric triple $(\fg,\theta,\ip_-)$ the extension $\ip$ of
$\ip_-$ in Condition (S3) is uniquely
determined. Therefore we will also call $(\fg,\theta,\ip)$ a symmetric triple.

A {\it metric Lie algebra} $(\fg,\ip)$ (or $\fg$ in abbreviated notation)
is a finite-dimensional Lie algebra $\fg$ together with a (non-degenerate) invariant
scalar product. A  metric Lie algebra with involution is a triple
$(\fg,\theta,\ip)$ (also abbreviated by $\fg$ or by $(\fg,\theta)$), where
$(\fg,\ip)$ is a metric Lie algebra and
$\theta$ is an involutive isometric automorphism. A metric Lie algebra  with
involution $(\fg,\theta,\ip)$ is
a symmetric triple if and only if (S1) or (S2) is satisfied.

Let $(M,g)$ be a pseudo-Riemannian symmetric space and $(\fg,\theta,
\ip)$ the associated symmetric triple. If $\fg=\fg_+\oplus\fg_-$ is
the corresponding decomposition of $\fg$, then $\fg_+$ is the Lie
algebra of the holonomy group of $(M,g)$ and the adjoint
representation of $\fg_+$ on $\fg_-$ is the holonomy representation. According
to de Rham's decomposition theorem a pseudo-Riemannian manifold is called
indecomposable, if its holonomy representation has no proper non-trivial
non-degenerate invariant subspace. Hence a simply-connected pseudo-Riemannian
symmetric space is indecomposable if and only if the symmetric triple associated
with this symmetric space is not the direct sum of two non-vanishing symmetric
triples. We will call such symmetric triples {\it indecomposable}. In particular
a simply-connected pseudo-Riemannian symmetric space is indecomposable if and
only if it is not the product of two non-trivial pseudo-Riemannian symmetric spaces.
However note that there are non simply-connected pseudo-Riemannian symmetric spaces 
which are not indecomposable in this sense, but which are not the product of two 
non-trivial pseudo-Riemannian symmetric spaces.

A pseudo-Riemannian symmetric space $(M,g)$ as well as its associated symmetric triple 
$(\fg,\theta,
\ip)$ are called
{\it semi-simple} if the Lie algebra $\fg$ is semi-simple.

Finally let us introduce the notion of a {\it homomorphism (isomorphism)} for
each of the above defined algebraic objects. An isomorphism of metric Lie
algebras is a Lie algebra isomorphism which is also an isometry with respect to
the given inner products. A homomorphism (isomorphism)
$\phi:(\fl_1,\theta_{\fl_1})\rightarrow(\fl_2,\theta_{\fl_2})$ of Lie algebras
with involution is a Lie algebra homomorphism (isomorphism) satisfying
$\phi\circ \theta_{\fl_1} =\theta_{\fl_2}\circ \phi$. An isomorphism of metric
Lie algebras with involution is an isomorphism of the underlying  metric Lie
algebras which is also an isomorphism of Lie algebras with involution.

\section{Quadratic cohomology of $(\fl,\thl)$-modules}\label{grishkov}

Let us first recall the notion of quadratic cohomology introduced in \cite{KO2}.
Let $\fl$ be a finite-dimensional Lie algebra. An {\it orthogonal $\fl$-module} is a
tuple $(\rho,\fa,\ip_{\fa})$ (also $\fa$ or $(\rho,\fa)$ in abbreviated notation) such 
that $\rho$ is a representation of
the Lie algebra
$\fl$ on the finite-dimensional real vector space $\fa$ and $\ip_\fa$ is
a scalar product on $\fa$ satisfying
$$\langle \rho(L)A_1,A_2\rangle_\fa + \langle A_1, \rho(L)A_2\rangle_\fa =0$$
for all $L\in\fl$ and $A_1,A_2\in \fa$.

For $\fl$ and (any $\fl$-module) $\fa$ we have the
standard cochain complex $(C^*(\fl,\fa),d)$, where $C^p(\fl,\fa)=
\Hom(\bigwedge ^p\fl,\fa)$ and we have the corresponding cohomology groups
$H^p(\fl,\fa)$.
Furthermore we have the standard cochain complex $(C^*(\fl),d)$, which arises from the 
one-dimensional trivial representation of $\fl$. Even if $\fa$ is one-dimensional we 
will distinguish between $C^*(\fl,\fa)$ and $C^*(\fl)$.

We have a product
$$C^p(\fl,\fa)\times C^q(\fl,\fa)\longrightarrow C^{p+q}(\fl)$$
defined by the composition
$$C^p(\fl,\fa)\times C^q(\fl,\fa)\stackrel{\wedge}{\longrightarrow}
C^{p+q}(\fl,\fa\otimes \fa) \stackrel{\ip_\fa}{\longrightarrow} C^p(\fl).$$

Let $p$ be even. Then the group of quadratic $(p-1)$-cochains is the group
$${\cal C} ^{p-1}_Q(\fl,\fa)=C^{p-1}(\fl,\fa)\oplus C^{2p-2}(\fl) $$
with  group operation defined by
$$ (\tau_1,\sigma_1)*(\tau_2,\sigma_2)=(\tau_1+\tau_2, \sigma_1
+\sigma_2 +\textstyle{\frac12} \langle \tau_1\wedge \tau_2\rangle)\,.$$
We consider now the set  
$${\cal Z} ^{p}_Q(\fl,\fa)=\{(\alpha,\gamma) \in C^{p}(\fl,\fa)\oplus
C^{2p-1}(\fl) \mid d\alpha=0,\
d\gamma=\textstyle{\frac12}\langle\alpha \wedge\alpha\rangle\} $$
whose elements are called quadratic $p$-cocycles. Then the group  
${\cal C} ^{p-1}_Q(\fl,\fa)$ acts on ${\cal Z} ^{p}_Q(\fl,\fa)$ by
$$(\alpha,\gamma)(\tau,\sigma)=\Big(\,\alpha +d\tau,\gamma +d\sigma
+\langle(\alpha +\textstyle{\frac12} d\tau)\wedge\tau\rangle\,\Big).$$
and we define the quadratic cohomology set ${\cal H}
^{p}_Q(\fl,\fa)~:={\cal Z} ^{p}_Q(\fl,\fa)/ {\cal C}
^{p-1}_Q(\fl,\fa)$. As usual, we denote the equivalence class of
$(\alpha,\gamma)\in {\cal Z} ^{p}_Q(\fl,\fa)$ in ${\cal H}
^{p}_Q(\fl,\fa)$ by $[\alpha,\gamma]$.

Now we consider pairs $(\fl_i,\fa_i)$, $i=1,2$, where $\fl_i$ are finite-dimensional Lie 
algebras and
$\fa_i=(\rho_i,\fa_i)$ are orthogonal $\fl_i$-modules. We say that
$(S,U):(\fl_{1},\fa_{1})\rightarrow (\fl_{2},\fa_{2})$ is a {\it
morphism of pairs} if $S:\fl_{1}\rightarrow \fl_{2}$ is a Lie algebra
homomorphism and $U:\fa_{2}\rightarrow\fa_{1}$ is an isometric
embedding such that
$$ U\circ \rho_2(S(L))=\rho_{1}(L)\circ U\,.$$
Note that $U$ maps in the reverse direction.

Now let $T:=(S,U):(\fl_{1},\fa_{1})\rightarrow (\fl_{2},\fa_{2})$ be a
morphism of pairs. For all $p\in\NN_{0}$ we define maps
$$\begin{array}{lll}
    T^*: C^{p}(\fl_{2})\longrightarrow C^{p}(\fl_{1}), &\qquad&
    T^*(\gamma):=S^{*}\gamma\\[1ex]
    T^*: C^{p}(\fl_{2},\fa_{2})\longrightarrow
    C^{p}(\fl_{1},\fa_{1}), &\qquad&
    T^*(\alpha):= U\circ (S^{*}\alpha)\,.
\end{array}  $$  
Then $T^*$ commutes with the differentials and
\begin{equation}\label{E2}
\langle T^*\alpha \wedge T^*\tau\rangle = T^* \langle \alpha\wedge\tau
\rangle
\end{equation}
holds for all $\alpha \in C^p(\fl,\fa)$ and  $\tau\in C^q(\fl,\fa)$.
In particular, ${\cal Z} ^{p}_Q(\fl,\fa)$ is invariant under $T^*\oplus T^*$.
Moreover,
$T^*\oplus T^*$ restricted to ${\cal C} ^{p-1}_Q(\fl,\fa)$ is a group
homomorphism, and
\begin{equation}
\label{E3}
(T^*\oplus T^*)(\,(\alpha,\gamma)(\tau,\sigma)\,)
=(T^*\alpha,T^*\gamma)(T^*\tau,T^*\sigma)
\end{equation}
holds for all $(\alpha,\gamma)\in{\cal Z} ^{p}_Q(\fl,\fa)$ and
$(\tau,\sigma)\in{\cal C} ^{p-1}_Q(\fl,\fa)$. Hence,
$$ T^*:{\cal H}^{p}_Q(\fl,\fa)\longrightarrow {\cal
H}^{p}_Q(\fl,\fa),\qquad T^*([\alpha,\gamma]):=[T^*\alpha,T^*\gamma]$$
is correctly defined.

Now let $(\fl,\theta_\fl)$ be a Lie algebra with involution.
\begin{de}
An orthogonal
$(\fl,\theta_\fl)$-module is a tuple
$(\rho,\fa,\ip_\fa,\theta_\fa)$, where
\begin{enumerate}
\item   $(\rho,\fa,\ip_\fa)$
is an orthogonal $\fl$-module for the Lie algebra $\fl$,
\item  $\theta_\fa$ is an
involutive isometry of $\fa$ compatible with $\rho$ and
$\theta_\fl$~:
\begin{equation}\label{E10}
 \tha \circ \rho(\theta_\fl(L))  =\rho(L)\circ \tha
\end{equation}
for all $L\in\fl$.
\end{enumerate}
Often we abbreviate the notation  by $\fa$ or $(\rho,\fa)$.
\end{de}

If $\fa$ is an orthogonal $(\fl,\theta_\fl)$-module, then $(\theta_\fl,\theta_\fa) 
:(\fl,\fa)\rightarrow(\fl,\fa)$ is a
morphism of pairs. As explained above $(\theta_\fl,\theta_\fa)$ defines involutions
\begin{equation} \label{E9}
\Theta:=(\theta_\fl,\theta_\fa)^*
\end{equation}
on $C^*(\fl,\fa)$, $C^*(\fl)$.
Let  $C^p(\fl,\fa)=C^p(\fl,\fa)_+ \oplus C^p(\fl,\fa)_-$ and $C^p(\fl)=C^p(\fl)_+ \oplus 
C^p(\fl)_-$be the eigenspace
decomposition with respect to $\Theta$. Obviously,
$${\cal C} ^{p-1}_Q(\fl,\fa)_+~:=C^{p-1}(\fl,\fa)_+\oplus
C^{2p-2}(\fl)_+ $$
is the space of fixed vectors of $\Theta\oplus\Theta: {\cal C}
^{p-1}_Q(\fl,\fa)\rightarrow {\cal C} ^{p-1}_Q(\fl,\fa)$. Moreover,
$${\cal Z} ^{p}_Q(\fl,\fa)_+~:={\cal Z} ^{p}_Q(\fl,\fa)\cap\,
(C^{p}(\fl,\fa)_+\oplus C^{2p-1}(\fl)_+)$$
is the set of fixed points of  $\Theta\oplus\Theta: {\cal Z}
^{p}_Q(\fl,\fa)\rightarrow {\cal Z} ^{p}_Q(\fl,\fa)$.
Using (\ref{E3}) we see that ${\cal C}^{p-1}_Q(\fl,\fa)_+$
acts on ${\cal Z} ^{p}_Q(\fl,\fa)_+$.
\begin{de}
Let $p$ be even. Then the set
$${\cal H}
^{p}_Q(\fl,\theta_{\fl},\fa):={\cal Z} ^{p}_Q(\fl,\fa)_+/ {\cal C}
^{p-1}_Q(\fl,\fa)_+.$$
is called quadratic cohomology of $(\fl,\thl)$ with values in the orthogonal
$(\fl,\thl)$-module $\fa$.
\end{de}
On the other hand, $\Theta$ acts on ${\cal H}^{p}_Q(\fl,\fa)$. The
next proposition compares the set of invariants of this action with
${\cal H}^{p}_Q(\fl,\theta_{\fl},\fa)$.
\begin{pr}\label{P31}
The map
\begin{eqnarray*}
{\cal H}^{p}_Q(\fl,\theta_{\fl},\fa)&\longrightarrow&
{\cal H}^{p}_Q(\fl,\fa)^\Theta:=\{\,[\alpha,\gamma]\in{\cal
  H}^{p}_Q(\fl,\fa) \mid \Theta([\alpha,\gamma])=[\alpha,\gamma]\,\}\\
\,[\alpha,\gamma] &\longmapsto & [\alpha,\gamma]
\end{eqnarray*}
is a bijection.
\end{pr}
\proof
Let $\proj_{+}$ and $\proj_{-}$ denote the projections with respect to
the decomposition $C^q(\fl,\fa)= C^q(\fl,\fa)_{+}\oplus C^q(\fl,\fa)_{-}$.
We abbreviate $\proj_+\beta=:\beta_+$, $\proj_-\beta=:\beta_-$ for
$\beta\in C^q(\fl,\fa)$. Then $d(\beta_{+})=(d\beta)_{+}$ and we write
only $d\beta_{+}$.

First we show that the map is injective. Suppose that
$(\alpha_1,\gamma_1), (\alpha_2,\gamma_2) \in {\cal Z} ^{p}_Q(\fl,\fa)_+$
and $[\alpha_1,\gamma_1]= [\alpha_2,\gamma_2] \in {\cal H}
^{p}_Q(\fl,\fa)$. Then we have $(\alpha_1,\gamma_1)(\tau,\sigma)=
(\alpha_2,\gamma_2)$ for a suitable element $(\tau,\sigma)\in{\cal
  C}^{p-1}_Q(\fl,\fa)$. Applying $\proj_+$ we obtain
\begin{eqnarray}
\alpha_2&=&\alpha_1+d\tau_{+} \label{E4}\\
\gamma_2&=&\gamma_1+ d\sigma_{+}+
\langle(\alpha_1+\textstyle{\frac12} d\tau)_{+}\wedge
\tau_{+}\rangle +\langle(\alpha_1+\textstyle{\frac12} d\tau)_{-}\wedge
\tau_{-}\rangle. \label{E5}
\end{eqnarray}
Since  $(\alpha_1,\gamma_1)(\tau,\sigma)=(\alpha_2,\gamma_2)$ implies
$\alpha_1+d\tau=\alpha_2$ we have $d\tau_{+}=d\tau$ and
$d\tau_{-}=0$.  Now (\ref{E4}) and (\ref{E5}) yield
$(\alpha_1,\gamma_1)(\tau_{+},\sigma_{+})=(\alpha_2,\gamma_2)$.
Since $(\tau_{+},\sigma_{+})\in {\cal C}^{p-1}_Q(\fl,\fa)_+$ we
obtain $[\alpha_1,\gamma_1]= [\alpha_2,\gamma_2] \in {\cal H}
^{p}_Q(\fl,\fa)_+$.

Now we prove that the map is surjective. Suppose
$[\Theta\alpha,\Theta\gamma]=[\alpha,\gamma] \in {\cal H}^{p}_Q(\fl,\fa)$. Then
there exists an element $(2\tau,2\sigma)\in {\cal
  C}^{p-1}_Q(\fl,\fa)$ such that
\begin{equation}\label{E6}
(\Theta\alpha,\Theta\gamma)=(\alpha,\gamma)(2\tau,2\sigma).
\end{equation}
Applying $\proj_+$ and $\proj_-$ to the first
component of (\ref{E6}) we obtain
$$d\tau_+=0,\ d\tau_-=-\alpha_-.$$
Therefore, applying $\proj_-$ to the second component of (\ref{E6})
gives
$$
-\gamma_-=\gamma_-+2d\sigma_-+\langle
(\alpha_-+d\tau_-)\wedge 2\tau_+\rangle+\langle(\alpha_++d\tau_+)\wedge
2\tau_-\rangle
= \gamma_- +2d\sigma_-+\langle \alpha_+\wedge 2\tau_-\rangle,
$$
hence
$$\gamma_-=-d\sigma_--\langle \alpha_+\wedge\tau_-\rangle.$$
Consequently,
\begin{eqnarray*}
(\alpha,\gamma)&=&(\alpha_++\alpha_-,\gamma_++\gamma_-)\\
&=&(\alpha_+-d\tau_-,
\gamma_+-d\sigma_--\langle\alpha_+\wedge\tau_-\rangle)\\
&=&(\alpha_+-d\tau_-,\gamma_+ -\textstyle{\frac12} d\tau_-\wedge
\tau_- -d\sigma_--\langle(\alpha_+
-\textstyle{\frac12}d\tau_-)\wedge\tau_-\rangle\,)\\
&=& (\alpha_+,\gamma_+ -\textstyle{\frac12} d\tau_-\wedge
\tau_-)(-\tau_-,-\sigma_-).
\end{eqnarray*}
We obtain $[\alpha,\gamma]=[\alpha_+,\gamma_+ -\textstyle{\frac12}
d\tau_-\wedge\tau_-]\in {\cal H}
^{p}_Q(\fl,\fa)$. This proves the assertion since $(\alpha_+,\gamma_+
-\textstyle{\frac12} d\tau_-\wedge\tau_-)\in {\cal Z} ^{p}_Q(\fl,\fa)_+$.
\qed

Let $(\fl_{i},\theta_{i})$, $i=1,2$, be Lie algebras with involution
and let $\fa_{i}$, $i=1,2$, be orthogonal $(\fl_{i},\theta_{\fl_{i}})$-modules. We will
say that
$(S,U): (\fl_{1},\theta_{\fl_{1}},\fa_{1})\rightarrow
(\fl_{2},\theta_{\fl_{2}},\fa_{2})$ is a {\it morphism of triples} if

\begin{enumerate}
\item $(S,U):(\fl_{1},\fa_{1})\rightarrow (\fl_{2},\fa_{2})$
is a morphism of pairs and
\item $S:\fl_{1}\rightarrow \fl_{2}$ and
$U:\fa_{2}\rightarrow \fa_{1}$ are homomorphisms of Lie algebras
with involution.
\end{enumerate}
If in addition $S$ and $U$ are isomorphisms, then $(S,U)$ is called {\it
isomorphism} of triples.
If $(S,U): (\fl_{1},\theta_{\fl_{1}},\fa_{1})\rightarrow
(\fl_{2},\theta_{\fl_{2}},\fa_{2})$ is a morphism of triples, then  
$(S,U)^*:C^{p}(\fl_{2},\fa_{2})\rightarrow
    C^{p}(\fl_{1},\fa_{1})$ induces a map
$$(S,U)^*: {\cal H}_{Q}^{p}(\fl_{2},\theta_{\fl_{2}},\fa_{2})\longrightarrow
{\cal H}_{Q}^{p}(\fl_{1},\theta_{\fl_{1}},\fa_{1}).$$

\section{Quadratic extensions}\label{leda}

\subsection{Definition}\label{S41}

Let $(\fl,\theta_\fl)$ be a Lie algebra with involution and let
$(\rho,\fa,\ip_\fa,\tha)$ be an orthogonal $(\fl,\theta_\fl)$-module.
\begin{de}
A quadratic extension of $(\fl,\thl)$ by $\fa$ is a tuple $(\fg,\theta,\fri,i,p)$,
where
\begin{enumerate}
\item $(\fg,\theta)$ is a metric Lie algebra with involution,
\item $\fri\subset\fg$ is a $\theta$-invariant isotropic ideal of $\fg$ and
\item $i$ and $p$ are homomorphisms of Lie algebras with involution constituting
an exact
  sequence
$$0\longrightarrow \fa \stackrel i{\longrightarrow} \fg/\fri \stackrel
p \longrightarrow \fl\longrightarrow 0\,,$$
such that
\begin{itemize}
\item[(i)] this exact sequence is consistent with the representation
  $\rho$ of $\fl$ on $\fa$ and
\item[(ii)] $i$ is an isometry from $\fa$ onto $\fri ^\perp/\fri$,
\end{itemize}
\end{enumerate}
\end{de}
\begin{re}\label{R41}
{\rm
\begin{enumerate} \item If $(\fg,\theta,\fri,i,p)$ is a quadratic extension of
$(\fl,\theta_\fl)$ by $\fa$, then  $(\fg,\fri,i,p)$ is a quadratic extension of
the Lie algebra $\fl$ by the orthogonal $\fl$-module $\fa$ in the sense of
\cite{KO2}.
\item Let $(\fg,\theta,\fri,i,p)$ be a quadratic extension of $(\fl,\thl)$ by
$\fa$. Let $\tilde p:\fg\rightarrow\fl$ be the composition
of the natural projection $\fg\rightarrow \fg/\fri$ with $p$.
Now let $p^*:=\fl^*\rightarrow\fg$ be the dual map of $\tilde p$, where we
identify $\fg^*$ with $\fg$ using the non-degenerate inner product on $\fg$.
This homomorphism is injective since
$\tilde p$ is surjective. Its image
equals $(\ker \tilde p)^\perp=\fri$.  
Hence $p^*$
determines a second exact sequence of Lie algebras with involution
$$
0\longrightarrow \fl^*{\longrightarrow}\fg
{\longrightarrow} \fg/\fri\longrightarrow 0\ ,
$$
where we consider $\fl^*$ as abelian Lie algebra with involution $\thl^*$.
In particular, $(\fg,\theta)$ can be considered as the result of two subsequent
extensions of Lie algebras with involution which satisfy certain compatibility 
conditions: first
we extend $(\fl,\thl)$ by $(\fa,\tha)$ and then we extend the
resulting Lie algebra by $(\fl^*,\thl^*)$.
\end{enumerate}
}
\end{re}
If $(\fg,\theta)$ is a metric Lie algebra with involution and  if
$\fri\in\fg$ is an isotropic $\theta$-invariant ideal such that
$\fri^\perp/\fri$ is abelian, then
the sequence
\begin{equation}\label{ifix}
0\longrightarrow \fri^\perp/\fri\stackrel{i}{\longrightarrow}\fg/\fri
\stackrel{p}{\longrightarrow} \fg/\fri^\perp\longrightarrow 0
\end{equation}
defines a quadratic extension of $\fg/\fri^\perp$ (with the induced involution)
by the orthogonal
module $\fri^\perp/\fri$. We call (\ref{ifix}) the canonical extension
associated with $(\fg,\theta,\fri)$.
\begin{de}
Two quadratic extensions $(\fg_j,\theta_j,\fri_j,i_j,p_j)$, $j=1,2$, of $(\fl,\thl)$
by $\fa$ are called equivalent if there exists an isomorphism
of metric Lie algebras with involution
$\Psi: (\fg_1,\theta_1) \rightarrow (\fg_2,\theta_2)$
which maps $\fri_1$ onto $\fri_2$ and satisfies
$$ \overline{\Psi}\circ i_1=i_2\qquad\mbox{ and }\qquad
p_2\circ\overline{\Psi}=p_1\ ,$$
where $\overline{\Psi}:\fg_1/\fri_1\rightarrow\fg_2/\fri_2$ is the
induced map.
\end{de}

\subsection{The standard model}\label{neptun}

\begin{de}\label{munkel}
  Let $(\fl,\thl)$ be a Lie algebra with involution and let
  $(\rho,\fa,\ip_\fa,\tha)$ be an orthogonal $(\fl,\thl)$-module. We
  consider the vector space $$\fd:=\fl^*\oplus\fa\oplus\fl$$ and define
  an inner product $\ip$ and an involutive endomorphism $\theta$ on $\fd$ by
\begin{eqnarray*}
 \langle Z_1+A_1+L_1,Z_2+A_2+L_2\rangle&:=& \langle A_1,A_2\rangle_\fa
+Z_1(L_2) +Z_2(L_1) \\
\theta(Z+A+L)&:=& \thl ^*(Z)+\tha(A)+\thl(L)
\end{eqnarray*}
for $Z,\,Z_1,\,Z_2\in \fl^*$, $A,\,A_1,\,A_2\in \fa$ and
$L,\,L_1,\,L_2\in \fl$. Now we choose $\alpha \in C^2(\fl,\fa)$ and
$\gamma \in C^3(\fl)$ and define an antisymmetric bilinear map
$\lb:\fd\times\fd\rightarrow \fd$ by
\begin{eqnarray*}
\ [\fl^*,\fl^*\oplus\fa] &=&0\\
\ [L_1,L_2] &=& \gamma(L_1,L_2,\cdot) +\alpha(L_1,L_2)+[L_1,L_2]_\fl\\
\ [L,A] &=& \rho(L)A - \langle A,\alpha(L,\cdot)\rangle\\
\ [L,Z]& = & \ad ^*(L)(Z)\\
\ [A_1,A_2]&=&\langle\rho(\cdot)A_1,A_2\rangle
\end{eqnarray*}
for $Z\in \fl^*$, $A,\,A_1,\,A_2\in \fa$ and
$L,\,L_1,\,L_2\in \fl$.
\end{de}
\begin{re}
{\rm
The triple $(\fd,\ip,\lb)$ coincides with the triple
$\fd_{\alpha,\gamma}(\fl,\fa,\rho)$ defined in
\cite{KO2}.
}
\end{re}

\begin{pr}
If $(\alpha,\gamma) \in {\cal Z}^2_Q(\fl,\fa)_+$, then
$\dd:=(\fd,\ip,\lb,\theta)$ is a metric Lie algebra with involution.
\end{pr}
\proof
This is a direct calculation. See \cite{KO2}, Section 3.2 for a proof
that  $(\alpha,\gamma) \in {\cal Z}_Q^2(\fl,\fa)$ implies that
$(\fd,\lb,\ip)$ is a metric Lie algebra. Obviously $\theta$ is an
isometry. It remains to prove that $\theta$ is a Lie algebra
homomorphism. Because of
\begin{eqnarray*}
 \lefteqn{ \theta ^{-1}[\theta(L_1),\theta(L_2)] =}\\
 &=& \thl^* (\,
  \gamma(\thl(L_1),\thl(L_2),\cdot)\,) +
 \tha(\,\alpha(\thl(L_1),\thl(L_2))\,)+\thl(\,[\thl(L_1),\thl(L_2)]_\fl\,)\\
&=& \gamma(\thl(L_1),\thl(L_2),\thl (\cdot)) +
 \tha(\,\alpha(\thl(L_1),\thl(L_2))\,)+[L_1,L_2]_\fl\\
&=&  (\Theta\gamma)(L_1,L_2,\cdot) + (\Theta\alpha)(L_1,L_2)+[L_1,L_2]_\fl\\
&=&  \gamma(L_1,L_2,\cdot) +\alpha(L_1,L_2)+[L_1,L_2]_\fl \ = \
 [L_1,L_2]
\end{eqnarray*}
we obtain $[\theta(L_1),\theta(L_2)]=\theta([L_1,L_2])$ for all
$L_1,\,L_2\in\fl$. The remaining identities
$[\theta(L),\theta(A)]=\theta([L,A])$,
$[\theta(L),\theta(Z)]=\theta([L,Z])$,
 and
$[\theta(A_1),\theta(A_2)]=\theta([A_1,A_2])$ for $Z\in \fl^*$, $A,A_1,\,A_2\in
\fa$ and
$L\in \fl$ can be proved in a
similar way using the compatibility of $\thl,\tha$ and $\rho$.
\qed

We identify $\fd/\fl^{*}$ with $\fa\oplus\fl$ and denote by
$i:\fa\rightarrow \fa\oplus\fl$ the injection and by $p:\fa\oplus\fl\rightarrow
\fl$ the projection. Then we have
\begin{co} If $(\alpha,\gamma) \in \cZ^{2}_{Q}(\fl,\fa)_+$, then
    $(\dd,\theta,\fl^{*},i,p)$ is a quadratic extension of
    $(\fl,\thl)$ by $(\fa,\tha)$.
\end{co}
We will denote the quadratic extension $(\dd,\theta,\fl^{*},i,p)$ also by
$\dd$.
\begin{re}\label{R51}{\rm
Let $\fa=\fa_+\oplus\fa_-$, $\fl=\fl_+\oplus\fl_-$, and $\fd=\fd_+\oplus\fd_-$ be the 
eigenspace decompositions with respect to the corresponding involutions.
If the signature of the restriction of $\ip_\fa$ to $\fa_-$ equals $(p,q)$, then the 
signature of $\ip$ restricted to $\fg_-$ equals $(p+\dim\fl_-,q+\dim\fl_-)$.
If $\dd$ is a symmetric triple, then this is also the signature of the metric on any 
pseudo-Riemannian symmetric space which is associated with this triple. }
\end{re}

\subsection{Classification by cohomology}

\begin{pr}\label{prop1}\quad
For $(\alpha_1,\gamma_1),\ (\alpha_2,\gamma_2) \in {\cal Z}^2_Q(\fl,\fa)_+$ the
quadratic extensions \linebreak $\fd_{\alpha_1,\gamma_1}(\fl,\thl,\fa)$ and
$\fd_{\alpha_2,\gamma_2}(\fl,\thl,\fa)$ of $(\fl,\thl)$ by $\fa$ are equivalent
if and only if $[\alpha_1,\gamma_1]=[\alpha_2,\gamma_2]\in{\cal
H}^2_Q(\fl,\thl,\fa)$.
\end{pr}
\proof
Assume first that $\Psi:\fd_{\alpha_1,\gamma_1}(\fl,\thl,\fa) \rightarrow
\fd_{\alpha_2,\gamma_2}(\fl,\thl,\fa)$ is an equivalence. The following facts
can be verified by direct calculations, see also \cite{KO2}, Prop.\,3.3 for a
detailed proof. Since $\Psi$ is an isometry and satisfies $\Psi(\fl^*)=\fl^*$,
$\proj_{\fa}\Psi|_{\fa}=\Id$ and
$\proj_{\fl}\Psi|_{\fl}=\Id$ it can be written as
\begin{equation}\label{EM1}
    \Psi=\left(
\begin{array}{ccc}
    \Id & -\tau^{*} & \bar\sigma -\frac 12\tau^{*}\tau  \\
    0 & \Id & \tau  \\
    0 & 0 & \Id
\end{array}\right) :\fl^{*}\oplus\fa\oplus\fl \longrightarrow
\fl^{*}\oplus\fa\oplus\fl,
\end{equation}
where $\tau\in C^{1}(\fl,\fa)$ and  $\sigma (\cdot\,,\cdot)= \langle\bar\sigma
(\cdot),\cdot\rangle \in C^{2}(\fl)$. Moreover, since $\Psi$ is a Lie algebra
isomorphism the cochains $\tau\in C^{1}(\fl,\fa)$ and  $\sigma\in C^{2}(\fl)$
satisfy $(\alpha_1,\gamma_1)(\tau,\sigma)=(\alpha_2,\gamma_2)$. Furthermore,
$\tau$ and $\sigma$ are $\Theta$-invariant, since $\Psi$ commutes with
$\theta=\thl^*\oplus\tha\oplus\thl$. Hence
$[\alpha_1,\gamma_1]=[\alpha_2,\gamma_2]\in{\cal H}^2_Q(\fl,\thl,\fa)$.
Conversely, if $[\alpha_1,\gamma_1]=[\alpha_2,\gamma_2]\in{\cal
H}^2_Q(\fl,\thl,\fa)$, then there exist cochains $\tau\in C^{1}(\fl,\fa)_+$ and
 $\sigma\in C^{2}(\fl)_+$ such that
$(\alpha_1,\gamma_1)(\tau,\sigma)=(\alpha_2,\gamma_2)$ holds and we can define
an equivalence $\Psi:\fd_{\alpha_1,\gamma_1}(\fl,\thl,\fa) \rightarrow
\fd_{\alpha_2,\gamma_2}(\fl,\thl,\fa)$ by (\ref{EM1}).
\qed

\begin{lm}\label{Ls}
Let $(\fg,\theta,\fri,i,p)$ be a quadratic extension of $(\fl,\theta_\fl)$ by
$\fa$ and let $\tilde p:\fg\rightarrow \fl$ be the composition of the natural
projection $\fg\rightarrow \fg/\fri$ with $p$. Then there exists an injective
homomorphism of vector spaces $s: \fl\rightarrow \fg$ such that
\begin{itemize}
\item[(i)] $\tilde p\circ s =\Id$,
\item[(ii)] $s\circ \thl =\theta\circ s$, and
\item[(iii)] $s(\fl)$ is isotropic.
\end{itemize}
\end{lm}
\proof Let us consider the orthogonal decomposition $\fg=\fg_+\oplus\fg_-$ into
eigenspaces of $\theta$. Since $\fri^\perp$ is $\theta$-invariant we have
$\fri^\perp =\fg_+\cap\fri^\perp\oplus \fg_-\cap\fri^\perp$. Moreover,
$\fri_+^\perp:=\fg_+\cap\fri^\perp$ and $\fri_-^\perp:=\fg_-\cap\fri^\perp$ are
coisotropic subspaces of $\fg_+$ and $\fg_-$, respectively. Therefore we can choose
isotropic vector space complements $V_+$ of $\fri_+^\perp$ in $\fg_+$ and $V_-$
of $\fri_-^\perp$ in $\fg_-$. Then $V:=V_+\oplus V_-$ is a $\theta$-invariant
isotropic vector space complement of $\fri^\perp$ in $\fg$. Hence, $s:=(\tilde 
p|_V)^{-1}:\fl\rightarrow V\subset \fg$ satisfies Conditions
$(i)$, $(ii)$, and $(iii)$ of the lemma.
\qed

\begin{pr}\label{prop2}\
Let $(\fg,\theta,\fri,i,p)$ be a quadratic extension of $(\fl,\theta_\fl)$ by
$(\rho,\fa)$ and let \linebreak $s: \fl\rightarrow \fg$ be as in Lemma \ref{Ls}. Let
$\alpha\in C^2(\fl,\fa)$ and $\gamma\in C^3(\fl)$ be defined by
\begin{eqnarray}
i(\alpha(L_1,L_2))&:=&[s(L_1),s(L_2)]-s([L_1,L_2])\, +\, \fri\quad
\in\,
\fg/\fri
\label{Ealpha}\\
\gamma(L_1,L_2,L_3)&:=& \langle\, [s(L_1),s(L_2)], s(L_3)\rangle \ .
\label{Egamma}
\end{eqnarray}
Then $(\alpha,\gamma)\in \cZ^2_Q(\fl,\fa)_+$ holds and the quadratic
extension  $(\fg,\theta,\fri,i,p)$ is equivalent to
$\dd$.
\end{pr}
\proof
By \cite{KO2}, Proposition 3.4 we already know that $(\alpha,\gamma)\in
\cZ^2_Q(\fl,\fa)$ and that $(\fg,\fri,i,p)$ and
the triple $\fd:=(\fd,\ip,\lb)$ associated with $(\alpha,\gamma)$ by Definition 
\ref{munkel} are equivalent as  quadratic extensions
 of the Lie algebra $\fl$ by the orthogonal $\fl$-module $\fa$ (disregarding 
involutions). In fact we proved the following. If
$p^*:\fl^*\rightarrow \fri$ is the isomorphism defined in Remark \ref{R41} and if we 
define
the linear map $t:\fa\rightarrow (\fri\oplus s(\fl))^\perp\subset \fg$  by
 $$i(A)=t(A)+\fri\in\fg/\fri,$$
then $$\Psi=p^* + t+s\ :\quad \fd=\fl^*\oplus \fa\oplus\fl\longrightarrow \fg,$$
 is an equivalence of $\fd$ and $(\fg,\fri,i,p)$.

It remains to show that $(\alpha,\gamma)\in \cZ^2_Q(\fl,\fa)$ is
$\Theta$-invariant and that $\Psi$ is compatible with the involutions
$\thl^*\oplus\tha\oplus\thl$ on $\fl^*\oplus\fa\oplus\fl$ and $\theta$ on $\fg$.
The first assertion follows easily from (\ref{Ealpha}) and (\ref{Egamma}) using
that $i$ is a homomorphism of Lie algebras with involution and that $s$
satisfies $s\circ \thl =\theta\circ s$. Let us now verify the second assertion.
Since $p$ is a homomorphism of Lie algebras with involution and $s$ satisfies
$s\circ \thl =\theta\circ s$ it remains to prove that $t\circ \tha =\theta\circ
t$ holds. Recall that $t\circ \tha(A)\in (\fri\oplus s(\fl))^\perp$ for all
$A\in\fa$. Since $\fri\oplus s(\fl)$ is $\theta$-invariant, also $\theta\circ
t(A)\in(\fri\oplus s(\fl))^\perp$ for all $A\in\fa$. Hence it suffices to prove
that $t\circ \tha(A)+\fri=\theta\circ t(A) +\fri\in\fg/\fri$ for all $A\in\fa$.
However, this is true by definition of $t$ and the fact that $i$
is a homomorphism of Lie algebras with involution.
 \qed
 
The following fact follows from Propositions \ref{prop1} and \ref{prop2}.
\begin{co} Let  $(\fg,\theta,\fri,i,p)$ be a quadratic extension of
$(\fl,\theta_\fl)$ by $\fa$ and let  $\alpha\in C^2(\fl,\fa)$ and $\gamma\in
C^3(\fl,\fa)$ be defined as in Proposition \ref{prop2}. Then the cohomology
class $[\alpha,\gamma]\in {\cal H}_Q^2(\fl,\thl,\fa)$ does not depend on the
choice of $s$.
\end{co}

Finally we obtain:
\begin{theo}\label{mark}
 The equivalence classes of quadratic extensions of a Lie algebra with
 involution $(\fl,\thl)$ by an orthogonal  $(\fl,\thl)$-module $\fa$
 are in one-to-one correspondence with elements of ${\cal H}_Q^2(\fl,\thl,\fa)$.
\end{theo}

\section{Admissible extensions}\label{poseidon}

In this section we will equip each symmetric triple $(\fg,\theta,\ip)$ without 
semi-simple ideals with the structure of a quadratic extension in a canonical way. 
Therefore we are particularly interested in those (so-called admissible) quadratic 
extensions which come from this canonical procedure. The main result is Theorem 
\ref{twain} which describes the subset ${\cal H}_Q^2(\fl,\thl,\fa)_\sharp\subset
{\cal H}_Q^2(\fl,\thl,\fa)$ corresponding to admissible quadratic extensions of a Lie 
algebra with involution
$(\fl,\theta_\fl)$ by an orthogonal module $\fa$.


First we recall the construction of the canonical isotropic ideal of a metric
Lie algebra $\fg$ from \cite{KO2}, Section 4.
Let $\fg$ be a finite-dimensional Lie algebra. Then there are chains of ideals
\begin{eqnarray*}
\{0\}=S_0(\fg)\subset S_1(\fg)\subset S_2(\fg)\subset&\dots&\subset
S_{l_+}(\fg)=\fg \\
\fg=R_0(\fg)\supset R_1(\fg)\supset R_2(\fg)\supset&\dots&\supset
R_{l_-}(\fg)=\{0\}
\end{eqnarray*}
which are defined inductively as follows: $S_k(\fg)$ is the largest
ideal of $\fg$ containing $S_{k-1}(\fg)$ such that the $\fg$-module $S_k(\fg)/ 
S_{k-1}(\fg)$ is
semi-simple.  $R_k(\fg)$ is the smallest ideal contained in $R_{k-1}(\fg)$ such that the 
$\fg$-module
$R_{k-1}(\fg)/ R_k(\fg)$ is semi-simple.
Then we set
\begin{equation}\label{herakles}
\fri(\fg):=\sum_{k=1}^{l_--1}S_k(\fg)\cap R_k(\fg)\ .
\end{equation}
If $\fri(\fg)$ is a metric Lie algebra, then $\fri(\fg)$ is isotropic and
$\fri(\fg)^\perp/ \fri(\fg)$ is abelian provided $\fg$ does not contain non-trivial
simple ideals (\cite{KO2}, Lemma 4.2). By construction, $\fri(\fg)$ is invariant under
all automorphisms of $\fg$, in particular under all involutions.

\begin{de}
A quadratic extension $(\fg,\theta,\fri,i,p)$ of $(\fl,\theta_\fl)$ by $\fa$ is
called balanced if $\fri=\fri(\fg)$. It is called admissible if it is balanced
and if $(\fg,\theta,\ip)$ is a symmetric triple.
\end{de}

\begin{pr}\label{schmidt}
    Any symmetric triple $(\fg,\theta,\ip)$ without semi-simple ideals has the
structure of
    a balanced, hence admissible,
    quadratic extension in a canonical way.
\end{pr}

\proof Take the canonical extension associated with $(\fg,\theta,\fri(\fg))$
defined by (\ref{ifix}).
\qed

\begin{de}\label{zwiesel}
Let $(\fl,\theta_\fl)$ be a Lie algebra with involution such that
\begin{itemize}
\item[$(T_1)$] $\quad [\fl_-,\fl_-]=\fl_+\ .$
\end{itemize}
Let $(\rho,\fa)$ be a semi-simple orthogonal $(\fl,\theta_\fl)$-module and let
$(\alpha,\gamma)\in
{\cal Z}_Q^2(\fl,\fa)_+$.
Then $\fa=\fa^\fl\oplus \rho(\fl)\fa$, and we have a
corresponding decomposition $\alpha=\alpha_0+\alpha_1$. We consider the following
conditions:
\begin{enumerate}
\item[$(A_0)$]
Let $L_0\in \fz(\fl)\cap \ker \rho$ be such that there exist
elements
$A_0\in \fa$ and $Z_0\in \fl^*$ satisfying
for all $L\in\fl$
\begin{enumerate}
\item[(i)] $\alpha(L,L_0)=\rho(L) A_0 $,
\item[(ii)] $\gamma(L,L_0,\cdot)=-\langle A_0,\alpha(L,\cdot)\rangle_\fa +\langle
Z_0, [L,\cdot]_\fl\rangle$ as an element of $\fl^*$,
\end{enumerate}
then $L_0=0$.
\item[$(B_0)$] The subspace $\alpha_0(\ker \lb_\fl)\subset
\fa^\fl$ is non-degenerate.
\item[$(T_{2})$] $\quad \fa^\fl_+=\alpha_0(\Ker \lb_{\fl_-})$.
\item[$(A_k)$] $(k\ge 1)$\\
Let $\fk\subset S(\fl)\cap R_k(\fl)$ be an $\fl$-ideal such that there exist
elements
$\Phi_1\in \Hom(\fk,\fa)$ and $\Phi_2\in \Hom(\fk,R_k(\fl)^*)$ satisfying
for all $L\in\fl$ and $K\in\fk$
\begin{enumerate}
\item[(i)] $\alpha(L,K)=\rho(L)\Phi_1(K)-\Phi_1([L,K]_\fl)$,
\item[(ii)] $\gamma(L,K,\cdot)=-\langle \Phi_1(K),\alpha(L,\cdot)\rangle_\fa
+\langle
\Phi_2(K), [L,\cdot]_\fl\rangle +\langle \Phi_2([L,K]_\fl), \cdot \rangle$ as an
element of $R_k(\fl)^*$,
\end{enumerate}
then $\fk=0$.
\item[$(B_k)$] $(k\ge 1)$\\
Let $\fb_k\subset\fa$ be the maximal $\fl$-submodule such that the system of equations
$$ \langle\alpha(L,K),
B\rangle_\fa=\langle\rho(L)\Phi(K)-\Phi([L,K]_\fl),B\rangle_\fa\ ,
\quad L\in\fl, K\in R_k(\fl), B\in\fb_k, $$
has a solution $\Phi \in \Hom(R_k(\fl),\fa)$. Then $\fb_k$ is non-degenerate.
\end{enumerate}
If the conditions
{$(T_{2})$}, $(A_k)$, $(B_k)$, $0\le k\le m$,
where $m$ is such that $R_{m+1}(\fl)=0$, are satisfied, then
the cohomology class $[\alpha,\gamma]\in {\cal H}_Q^2(\fl,\thl,\fa)$ is called 
admissible. We denote the set of all
admissible cohomology classes by ${\cal H}_Q^2(\fl,\thl,\fa)_\sharp$. A
Lie algebra
with involution $(\fl,\theta_\fl)$ is called admissible if it satisfies $(T_1)$ and there is a
semi-simple
orthogonal $(\fl,\theta_\fl)$-module $\fa$ such that
${\cal H}_Q^2(\fl,\thl,\fa)_\sharp\ne\emptyset$.
\end{de}

Each of the above conditions depends only on the cohomology class $[\alpha,\gamma]$
and not on the particular cocycle $(\alpha,\gamma)$ representing it (compare the
discussion
in \cite{KO2}, Section~4). Thus the notion of admissibility of a cohomology
class is well-defined. In addition, we remark that the submodule $\fb_k$
appearing in Condition $(B_k)$ is automatically $\tha$-invariant.

In the remainder of this section we want to prove the following theorem.

\begin{theo}\label{twain}
Let $(\fl,\theta_\fl)$ be a Lie algebra with involution.
Let $(\rho,\fa)$ be an orthogonal $(\fl,\theta_\fl)$-module and let
$[\alpha,\gamma]\in
{\cal H}_Q^2(\fl,\thl,\fa)$. Then the quadratic extension $\dd$ is admissible if
and only
if
\begin{enumerate}
\item[(i)] $(\fl,\theta_\fl)$ satisfies $(T_1)$,
\item[(ii)] the representation $\rho$ is semi-simple, and
\item[(iii)] $[\alpha,\gamma]\in
{\cal H}_Q^2(\fl,\thl,\fa)_\sharp$.
\end{enumerate}
\end{theo}

In the following we abbreviate $\fd:=\dd$, $\fd=\fd_+\oplus\fd_-$.
\begin{lm}\label{p52}
Let $\rho$ be semi-simple. If $[\fd_-,\fd_-]=\fd_+$, then Conditions
$(T_1)$ and $(T_{2})$ are satisfied.
\end{lm}
\proof The only non-trivial assertion to be proved is $\fa^\fl_+\subset\alpha_0(\Ker
\lb_{\fl_-})$. Suppose $A\in\fa_+^\fl$. By assumption there are elements
$X_i^1=Z^1_i+A_i^1+L_i^1$, $X_i^2=Z^2_i+A_i^2+L_i^2$,  $i=1,...,k$, in $\fd_-$
such that $$A=\sum_{i=1}^k[X_i^1,X_i^2] = \sum_{i=1}^k[L_i^1,L_i^2]_\fl
+\sum_{i=1}^k
\Big(\alpha(L_i^1,L_i^2)+\rho(L_i^1)(A_i^2)-\rho(L_i^2)(A_i^1)\Big)+Z $$
with $Z\in\fl^*$. Obviously, $Z=0$. Moreover,
$\rho(L_i^1)(A_i^2)-\rho(L_i^2)(A_i^1)=0$ since this term is in $\rho(\fl)\fa$
and $\fa=\fa^\fl\oplus \rho(\fl)\fa$. Hence,
$A= \sum_{i=1}^k \alpha(L_i^1,L_i^2)$ and $\sum_{i=1}^k[L_i^1,L_i^2]_\fl=0$.
\qed
\begin{lm} \label{P52}
Let $\rho$ be semi-simple. If the Conditions $(T_1)$, $(T_2)$, {$(A_0)$}
are satisfied, then
the representation of $\fd_+$ on $\fd_-$ is faithful.
\end{lm}
\proof First we note that $(T_1)$ implies
\begin{equation}
\label{UL1}
(\rho(\fl)\fa)_+=\rho(\fl_-)\fa_-.
\end{equation}
To verify this we have to prove that $L(A)\in \rho(\fl_-)\fa_-$ for all
$L\in\fl_+$ and $A\in\fa_+$. Because of  $(T_{1})$ it suffices to show that
$[L_1,L_2](A)\in \rho(\fl_-)\fa_-$ for all $L_1,L_2\in\fl_-$ and $A\in\fa_+$.
But this is obvious since $[L_1,L_2](A)=L_1L_2(A)-L_2L_1(A)$ and
$L_i(A)\in\fa_-$ for $i=1,2$.

Equation (\ref{UL1}) implies
$$(\rho(\fl)\fa)_+ +\fl^*_+ = [\fl_-,\fa_-]+\fl^*_+ \subset
[\fd_-,\fd_-]+\fl^*_+, $$
whereas $(T_{2})$ yields
$$\fa_+^\fl \subset[\fl_-,\fl_-]+\fl^*_+ \subset [\fd_-,\fd_-]+\fl^*_+.$$
Consequently,
$$\fa_+\subset[\fd_-,\fd_-]+\fl^*_+$$ and therefore
\begin{equation}\label{U1}
\fd_+=[\fl_-,\fl_-]+\fa_++\fl_+^* = [\fd_-,\fd_-]+\fl_+^*.
\end{equation}

Assume now, that $Z_+\in \fl^*_+$, $A_+\in\fa_+$ and $L_+\in\fl_+$ are such that
$[Z_++A_++L_+,\fd_-]=0$. In particular this implies
$[L_+,\fl_-]_\fl=0$. Because of $(T_{1})$ this yields $L_+\in\fz(\fl)$. Using
this and (\ref{U1}) we obtain
$$[Z_++A_++L_+,\fd_+]= [Z_++A_++L_+,\fl_+^*]=\ad^*(L_+)(\fl_+^*)=0,$$
hence $Z_++A_++L_+\in\fz(\fd)$.
A straightforward computation shows that Condition $(A_0)$ is equivalent to
$\fz(\fd)\subset
\fl^*\oplus\fa$. It follows that $L_+=0$. Hence we have $Z_++A_+\in\fz(\fd)$. In
particular, $[Z_++A_+,\fl]=0$ holds. This implies $A_+\in\fa^\fl$ and
$$\langle A_+,\alpha(\cdot\,,\cdot)\rangle =\langle
A_+,\alpha_0(\cdot\,,\cdot)\rangle =-Z_+(\lb).$$
However the last equation yields $A_+\perp \alpha_0(\Ker \lb_{\fl_-})=\fa_+^\fl$, hence
$A_+=0$. We obtain $Z_+\in\fz(\fd)$ and, in particular, $[Z_+,\fl]=0$, which
yields $Z_+([\fl,\fl]_\fl)=0$. This implies $Z_+(\fl_+)=0$, thus $Z_+=0$.
Consequently, the representation of $\fd_+$ on $\fd_-$ is faithful.
\qed

{\sl Proof of Theorem \ref{twain}. }
Theorem 4.1 in \cite{KO2} tells us that the quadratic extension $\dd$ is
balanced if and only
if $(\rho,\fa)$ is semi-simple and the conditions $(A_k)$, $(B_k)$, $0\le k\le
m$, are satisfied.
Recall from Section \ref{sym} that a metric Lie algebra with involution is a
symmetric triple
if and only if one of the Conditions {\rm(S1)} and {\rm (S2)} is satisfied. The
theorem now follows from Lemma \ref{p52} and Lemma \ref{P52}.
\qed

\section{A classification scheme for indecomposable symmetric triples}\label{quark}

By the results of the previous section the metric Lie algebras
with involution $\dd$
associated
with semi-simple orthogonal $(\fl,\thl)$-modules $\fa$ of an admissible Lie
algebra with involution $(\fl,\thl)$ and
$[\alpha,\gamma]\in \cH^2_Q(\fl,\thl,\fa)_\sharp$ exhaust all isomorphism classes
of symmetric triples without semi-simple ideals. It remains to decide
which of these data lead to isomorphic symmetric triples.
Since one is interested in the classification of {\em indecomposable}
symmetric triples we also have to check indecomposability of $\dd$ as
a symmetric triple in terms of the defining data $(\fl,\thl,\fa)$,
$[\alpha,\gamma]\in \cH^2_Q(\fl,\thl,\fa)_\sharp$. A comprehensive
discussion of the analogous questions for metric Lie algebras can
be found in \cite{KO2}, Section 5. This allows us to be rather brief
here. We will conclude the section with a classification scheme for
indecomposable symmetric triples.

\begin{pr}\label{wiesel}
Let $(\fl_i,\theta_{\fl_i})$, $i=1,2$, be
admissible Lie algebras with involution. Let $\fa_i$, $i=1,2$, be orthogonal
$(\fl_i,\theta_{\fl_i})$-modules
and suppose $(\alpha_i,\gamma_i)\in \cZ^2_Q(\fl_i,\fa_i)_+$ such that
$[\alpha_i,\gamma_i]\in \cH^2_Q(\fl_i,\theta_{\fl_i},\fa_i)_\sharp$.

Then $\fd_{\alpha_1,\gamma_1}(\fl_1,\theta_{\fl_1},\fa_1)$ and
$\fd_{\alpha_2,\gamma_2}(\fl_2,\theta_{\fl_2},\fa_2)$ are isomorphic as
symmetric triples,
if and only if there
is an isomorphism of triples $(S,U): (\fl_1,\theta_{\fl_1},\fa_1)\rightarrow
(\fl_2,\theta_{\fl_2},\fa_2)$ such that
$ (S,U)^*[\alpha_2,\gamma_2]=[\alpha_1,\gamma_1]\in
\cH^2_Q(\fl_1,\theta_{\fl_1},\fa_1)_\sharp$.
\end{pr}
\proof
If $(S,U): (\fl_1,\theta_{\fl_1},\fa_1)\rightarrow
(\fl_2,\theta_{\fl_2},\fa_2)$ is an isomorphism of triples, then
the symmetric triple $\fd_{\alpha_2,\gamma_2}(\fl_2,\theta_{\fl_2},\fa_2)$ is
isomorphic to $\fd_{U\circ(S^*\alpha_2), S^*\gamma_2}(\fl_1,\theta_{\fl_1},\fa_1)$.
If in addition $ (S,U)^*[\alpha_2,\gamma_2]=[\alpha_1,\gamma_1]$, then
the latter quadratic extension is equivalent to
$\fd_{\alpha_1,\gamma_1}(\fl_1,\theta_{\fl_1},\fa_1)$ by Proposition
\ref{prop1}. Thus $\fd_{\alpha_1,\gamma_1}(\fl_1,\theta_{\fl_1},\fa_1)$ and
$\fd_{\alpha_2,\gamma_2}(\fl_2,\theta_{\fl_2},\fa_2)$ are isomorphic.

For the reverse direction we really need the admissibility assumption.
Since both qua\-dratic extensions are balanced, an isomorphism
$\Psi:\fd_{\alpha_1,\gamma_1}(\fl_1,\theta_{\fl_1},\fa_1)\rightarrow
\fd_{\alpha_2,\gamma_2}(\fl_2,\theta_{\fl_2},\fa_2)$ maps the canonical
isotropic ideals into each other, hence is compatible
with the filtrations given by $\fl_i^*\subset \fl_i^*\oplus \fa_i$, $i=1,2$. It
therefore induces maps $S:\fl_1\rightarrow\fl_2$,
$U^{-1}:\fa_1\rightarrow \fa_2$. Then $(S,U)$ is a morphism of triples
and the quadratic extensions $\fd_{U\circ (S^*\alpha_2),
S^*\gamma_2}(\fl_1,\theta_{\fl_1},\fa_1)$ and
$\fd_{\alpha_1,\gamma_1}(\fl_1,\theta_{\fl_1},\fa_1)$ are equivalent.
Now we apply Proposition~\ref{prop1}.

For more details we refer to \cite{KO2}, Section 5.
\qed

\begin{de}\label{kosel}
A non-trivial decomposition of a triple $(\fl,\thl,\fa)$
consists of two non-zero morphisms of triples
$$ (q_i,j_i):(\fl,\thl,\fa)\longrightarrow (\fl_i,\theta_{\fl_i},\fa_i)\ ,
\quad i=1,2 \ ,$$
such that $(q_1,j_1)\oplus (q_2,j_2):(\fl,\thl,\fa)\rightarrow
(\fl_1,\theta_{\fl_1},\fa_1)\oplus(\fl_2,\theta_{\fl_2},\fa_2)$ is
an isomorphism.

A cohomology class $\varphi\in
\cH^2_Q(\fl,\theta_{\fl},\fa)$ is called decomposable if it can
be written as a sum
$$ \varphi=(q_1,j_1)^*\varphi_1 +(q_2,j_2)^*\varphi_2 $$
for a non-trivial decomposition $(q_i,j_i)$ of $(\fl,\thl,\fa)$ and
$\varphi_i\in \cH^2_Q(\fl_i,\theta_{\fl_i},\fa_i)$, $i=1,2$.
Here addition is induced by addition in the vector space $C^{p}(\fl,\fa)_+\oplus
C^{2p-1}(\fl)_+$.

A cohomology class which is not decomposable is called indecomposable.
We denote the set of all indecomposable elements in
$\cH^2_Q(\fl,\theta_{\fl},\fa)_\sharp$ by
$\cH^2_Q(\fl,\theta_{\fl},\fa)_0$.
\end{de}

\begin{pr}\label{hermlin}
Let $(\fl,\theta_{\fl})$ be an
admissible Lie algebra with involution. Let $\fa$ be an orthogonal
$(\fl,\theta_{\fl})$-module
and let $(\alpha,\gamma)\in \cZ^2_Q(\fl,\fa)_+$ be such that
$[\alpha,\gamma]\in \cH^2_Q(\fl,\theta_{\fl},\fa)_\sharp$.

Then $\fd_{\alpha,\gamma}(\fl,\theta_{\fl},\fa)$ is indecomposable if and
only if $[\alpha,\gamma]\in \cH^2_Q(\fl,\theta_{\fl},\fa)_0$.
\end{pr}
\proof
Assume that $[\alpha,\gamma]=(q_1,j_1)^*[\alpha_1,\gamma_1]
+(q_2,j_2)^*[\alpha_2,\gamma_2]$ for some non-trivial decomposition
of $(\fl,\thl,\fa)$. Note that
$$(q_1,j_1)^*[\alpha_1,\gamma_1]
+(q_2,j_2)^*[\alpha_2,\gamma_2]=((q_1,j_1)\oplus
(q_2,j_2))^*[(\alpha_1,\alpha_2),(\gamma_1,\gamma_2)]\ .$$
Arguing  as in the
first part of the proof of Proposition \ref{wiesel} we find that
$\dd$ is isomorphic to the orthogonal direct sum
$\fd_{\alpha_1,\gamma_1}(\fl_1,\theta_{\fl_1},\fa_1)\oplus
\fd_{\alpha_2,\gamma_2}(\fl_2,\theta_{\fl_2},\fa_2)$ and thus decomposable.

It is again only the reverse direction where we need the admissibility
assumptions. Let $\dd$ be decomposable.
Thus it is isomorphic to an orthogonal direct sum $\fd_1\oplus\fd_2$.
Without loss of generality we may assume
$\fd_i=\fd_{\alpha_i,\gamma_i}(\fl_i,\theta_{\fl_i},\fa_i)$
for certain triples $(\fl_i,\theta_{\fl_i},\fa_i)$ and
$[\alpha_i,\gamma_i]\in\cH^2_Q(\fl_i,\theta_{\fl_i},\fa_i)_\sharp$.
Then
$$[(\alpha_1,\alpha_2),(\gamma_1,\gamma_2)]
\in\cH^2_Q(\fl_1\oplus\fl_2,\theta_{\fl_1}\oplus\theta_{\fl_2},\fa_1\oplus\fa_2)_\sharp\
.$$
Now Proposition \ref{wiesel} implies that there exists an
isomorphism of triples
$$(S,U): (\fl,\thl,\fa)\rightarrow
(\fl_1\oplus\fl_2,\theta_{\fl_1}\oplus\theta_{\fl_2},\fa_1\oplus\fa_2)$$
such that
$(S,U)^*[(\alpha_1,\alpha_2),(\gamma_1,\gamma_2)]
=[\alpha,\gamma]$. Setting $q_i:=\proj_{\fa_i}\circ S$, $j_i:=U|_{\fa_i}$
we obtain the decomposition $[\alpha,\gamma]=(q_1,j_1)^*[\alpha_1,\gamma_1]
+(q_2,j_2)^*[\alpha_2,\gamma_2]$. Thus $[\alpha,\gamma]$ is decomposable.
\qed

Let us fix an admissible Lie algebra with involution $(\fl,\thl)$ and a
pseudo-Euclidean vector space $(\fa,\tha)$ with involution.
Conjugation by $\tha$ defines involutions on the Lie algebra
$\so(\fa)$ and on the group $O(\fa)$, also denoted by $\tha$.
We set
$$\Aut(\fa):=O(\fa)^{\tha}=O(\fa_+)\times O(\fa_-)\ .$$
We consider the set $\Hom((\fl,\thl),\so(\fa))_{\rm ss}$ of all orthogonal
semi-simple representations
of $\fl$ on $\fa$ which are compatible with the involutions $\thl$ and $\tha$.
If $\rho\in \Hom((\fl,\thl),\so(\fa))_{\rm ss}$ is fixed
we denote the corresponding orthogonal $(\fl,\thl)$-module by $\fa_\rho$.
The group $G:=\Aut(\fl,\thl)\times \Aut(\fa)$ acts from the right on
$\Hom((\fl,\thl),\so(\fa))_{\rm ss}$
by
$$(S,U)^*\rho:=\Ad(U^{-1})\circ S^*\rho\ ,\qquad S\in \Aut(\fl,\thl),\  U\in
\Aut(\fa) \ .$$
Then for any $\rho\in \Hom((\fl,\thl),\so(\fa))_{\rm ss}$ an element
$g=(S,U)\in G$ defines an isomorphism of triples
$ \bar g:=(S,U^{-1}): (\fl,\thl,\fa_{g^*\rho})\rightarrow (\fl,\thl,\fa_\rho)$
and therefore induces a bijection
$$ \bar g^*: \cH^2_Q(\fl,\thl,\fa_\rho)\rightarrow \cH^2_Q(\fl,\thl,\fa_{g^*\rho})\
.$$
We
obtain a right action of $G$ on the disjoint union
$$ \coprod_{\rho\in \Hom((\fl,\thl),\so(\fa))_{\rm ss}}
\cH^2_Q(\fl,\thl,\fa_\rho)\ .$$
As in Definition \ref{kosel} let $\cH^2_{Q}(\fl,\thl,\fa_\rho)_{0}\subset
\cH^2_Q(\fl,\fa_\rho)$ be the subset
of all admissible indecomposable elements (see Definition \ref{zwiesel} for the
admissibility conditions). Then the set
$$ \coprod_{\rho\in \Hom((\fl,\thl),\so(\fa))_{\rm ss}}
\cH^2_Q(\fl,\thl,\fa_\rho)_{0}$$
is $G$-invariant.
Combining Proposition \ref{schmidt} and Theorem \ref{twain} with Propositions
\ref{wiesel} and \ref{hermlin} we obtain
\begin{theo}\label{Pwieder}
Let $(\fl,\thl)$ be an admissible Lie algebra with involution, and
let $(\fa,\tha)$ be a pseudo-Euclidean vector space with involution. We consider the
class $\cA(\fl,\thl,\fa)$ of non-semisimple indecomposable symmetric
triples $(\fg,\theta,\ip)$
                       satisfying
    \begin{enumerate}
        \item The Lie algebras with involution $\fg/\fj(\fg)$ and
        $(\fl,\thl)$ are isomorphic.
        \item $\fj(\fg)/\fri(\fg)$ is isomorphic to $(\fa,\tha)$
	as a pseudo-Euclidean vector space with involution.
    \end{enumerate}
    Then the set of isomorphism classes of $\cA(\fl,\thl,\fa)$ is in bijective
    correspondence with the orbit space of the action of $G=\Aut(\fl,\thl)\times
    \Aut(\fa)$ on
    $$\coprod_{\rho\in \Hom((\fl,\thl),\so(\fa))_{\rm ss}}
\cH^2_Q(\fl,\thl,\fa_\rho)_{0}\
.$$
This orbit space can also be written as
$$\coprod_{[\rho]\in \Hom((\fl,\thl),\so(\fa))_{\rm ss}/G}
\cH^2_Q(\fl,\thl,\fa_\rho)_{0}/G_\rho\ ,$$
where $G_\rho=\{g\in G\:|\: g^*\rho=\rho\}$ is the automorphism group of the triple
$(\fl,\thl,\fa_\rho)$.
\end{theo}

\section{Classification of symmetric triples of index 2}
\label{S7}

Our general classification principle yields more explicit results if one only considers 
symmetric triples of a (fixed) small index.
Here the {\it index of a symmetric triple} $(\fg,\theta,\ip_-)$ is simply the index of 
$\ip_-$ or, equivalently, the index of the metric of any pseudo-Riemannian symmetric 
space associated with $(\fg,\theta,\ip_-)$.

In the following section we will demonstrate this for symmetric triples of index 2. We 
will obtain a list of all isomorphism classes of indecomposable symmetric triples of 
index 2 or, equivalently, a list of all isometry classes of indecomposable 
simply-connected pseudo-Riemannian symmetric spaces of index 2. First classification 
results for this case were already obtained by M.\,Cahen and M.\,Parker in \cite{CP1} 
and \cite{CP}. In \cite{CP1} symmetric spaces of index 2 with solvable transvection 
group were studied. Unfortunately, the classification of these spaces was not complete. 
In his diploma thesis Th.\,Neukirchner elaborated the claimed results, found the gaps 
and gave a revised classification of symmetric spaces of index 2 with solvable 
transvection group. Comparing the latter one  with the classification which follows from 
our classification scheme we observed that also Neukirchner's classification is not 
quite correct
(besides minor errors a series of spaces is missing and some of the normal forms contain 
too much parameters).

Now let $(\fg,\theta,\ip)$ be an indecomposable symmetric triple of index~2 which is not 
semi-simple. By Prop.\ \ref{schmidt} $(\fg,\theta,\ip)$ has the structure of a quadratic 
extension $\dd$ in a canonical way, where $(\fl,\thl)$ is a suitable admissible Lie 
algebra with involution, $\fa$ is a suitable orthogonal $(\fl,\thl)$-module and 
$[\alpha,\gamma]\in\cH^2_Q(\fl,\thl,\fa)_{0}$.
Using Remark \ref{R51} we see that $\dim\fl_-\le 2$. Since $(\fl,\theta_\fl)$ is 
admissible we have $\fl=[\fl_-,\fl_-]\oplus \fl_-$ and therefore $\dim \fl\le 3$.

In particular it follows from dimensional reasons that either $\fl$ is solvable or 
isomorphic to $\fsu(2)$ or to $\fsl(2,\RR)$. If $\fl$ is solvable, then also $\fg$ is 
solvable and we are in the case discussed above. If $\fl$ is isomorphic to $\fsu(2)$ or 
to $\fsl(2,\RR)$, then $\fl$ is the Levi factor of $\fg$. Symmetric triples with a Levi 
factor of this type were thoroughly studied in \cite{CP}. In particular, Cahen and 
Parker obtained a
classification of those triples
whose Lie algebra structure satisfies in addition a certain minimality condition.
In the special case of index~2 this minimality condition is always satisfied.
For completeness we will reproduce also this classification of non-solvable 
indecomposable
symmetric triples of index 2 using our classification scheme (taking only
the rather elementary classification of indecomposable orthogonal $(\fl,\thl)$-modules
from \cite{CP}).


\subsection{Semi-simple orthogonal representations of solvable Lie algebras with 
involution}

Before we concentrate on symmetric triples of index 2 we will give the following 
description of $\Hom((\fl,\thl),\so(\fa))_{\mbox{\footnotesize{\rm ss}}}/
\Aut (\fa)$ for an arbitrary solvable Lie algebra with involution
$(\fl,\thl)$ satisfying $[\fl_-,\fl_-]=\fl_+$ and the standard pseudo-Euclidean space
$\fa=\RR^{p+r,q+s}$ with involution
$\tha$ given by $\fa_+=\RR^{p,q}$, $\fa_-=\RR^{r,s}$. Here we have
$\Aut (\fa)=O(p,q)\times O(r,s)$.
Since $\fl'=[\fl,\fl]=R(\fl)$ and $R(\fl)$ is in the kernel of any
semi-simple representation of $\fl$ the map
\begin{eqnarray*}
\Hom((\fl/\fl',\bar{\theta}_\fl),\so(\fa))_{\mbox{\footnotesize{\rm ss}}}
&\longrightarrow&
\Hom((\fl,\thl),\so(\fa))_{\mbox{\footnotesize{\rm ss}}}\\
\bar\rho&\longmapsto&\rho,\ \rho(L)=\bar\rho(L+\fl)
\end{eqnarray*}
is a bijection. Here $\bar\theta_\fl$ is the involution induced by $\thl$ on the
quotient. Since by assumption $\fl_+\subset \fl'$ we see that
$\bar\theta_\fl=-\Id$. Hence for our purpose it is sufficient to determine
$\Hom((\fl,-\Id),\so(\fa))_{\mbox{\footnotesize{\rm ss}}}/
\Aut (\fa)$ for any abelian Lie algebra $\fl$.

Let $\RR^{p,q}$ be the standard pseudo-Euclidean space of dimension $n=p+q$.
Then the standard basis $e_1,\dots,e_n$ of $\RR^{p,q}$ is orthogonal and
satisfies $\langle e_k, e_k\rangle =-1$ for $k=1,\dots,p$ and $\langle
e_k,e_k\rangle=1$ for $k=p+1,\dots,n$.

We consider $\fa=\RR^{p+r,q+s}$ with the involutive isometry $\tha$ defined above by
$\fa=\fa_+\oplus\fa_-$, where $\fa_+ =\RR^{p,q}$ and  
$\fa_- =\RR^{r,s}$. We will use the notation $\fa_+^{p,q}$ for $\fa_+$, $\fa_-^{r,s}$ 
for $\fa_-$ and
$\fa_+^{p,q}\oplus\fa_-^{r,s}$ for the pair $(\fa,\tha)$. We will say that 
$A_1,\dots,A_{p+q+r+s}$ is a standard basis of $\fa_+^{p,q}\oplus\fa_-^{r,s}$ if 
$A_1,\dots A_{p+q}$ is the standard basis of $\fa_+^{p,q}=\RR^{p,q}$ and 
$A_{p+q+1},\dots A_{p+q+r+s}$ is the standard basis of $\fa_-^{r,s}=\RR^{r,s}$. 
given
$n=p+q+r+s$, $n'=p'+q'+r'+s'$, is a standard basis of  
if $A_1,\dots, A_n$ is the standard basis of $\fa_+^{p,q}\oplus\fa_-^{r,s}$ and 
$A_{n+1},\dots, A_{n+n'}$ is the standard 
Often we identify
$(\fa\oplus\fa',\tha\oplus\tha')$ with $\fa_+^{p+p',q+q'}\oplus\fa_-^{r+r',s+s'}$.

Now let $(\fl,\thl)$ be an abelian Lie algebra with involution $\thl=-\Id$. For
$\lambda\in\fl^*$ we define orthogonal representations $\rho_\lambda^+$ of
$(\fl,\thl)$ on $\fa_+^{0,1}\oplus\fa_-^{0,1}$,  
$\rho_\lambda^-$ of $(\fl,\thl)$ on $\fa_+^{1,0}\oplus\fa_-^{1,0}$,
$\tilde\rho^+_\lambda$ of $(\fl,\thl)$ on $\fa_+^{1,0}\oplus\fa_-^{0,1}$, and
$\tilde\rho^-_\lambda$ of $(\fl,\thl)$ on $\fa_+^{0,1}\oplus\fa_-^{1,0}$ by
\begin{eqnarray*}
&\rho^\pm_{\lambda}(L)=\small{
\left(\begin{array}{cc}
   0&-\lambda(L)\\
   \lambda(L)&0
   \end{array}
\right)},\quad
\tilde\rho^\pm_{\lambda}(L)=\small{
\left(\begin{array}{cc}
   0&\lambda(L)\\
   \lambda(L)&0
   \end{array}
\right)}&
\end{eqnarray*}
w.\,r.\,t.\ the standard bases.
Moreover, for $\mu,\,\nu\in\fl^{*}$ we define an orthogonal representation
$\rho''_{\mu,\nu}$ of
$(\fl,\thl)$ on $\fa_+^{1,1}\oplus\fa_-^{1,1}$ by
$$
\rho''_{\mu,\nu}(L)=\small{
        \left(
        \begin{array}{cccc}
            0 & 0 & -\nu(L) & \mu(L)  \\
            0 & 0 &\mu(L) &  \nu(L) \\
            \nu(L)  & \mu (L) & 0 & 0\\
            \mu(L)  & -\nu(L) & 0 & 0
        \end{array}
            \right)}
$$
w.\,r.\,t.\ the standard basis.

For
$\lambda=(\lambda^{1},\ldots,\lambda^{m}),\,\mu=(\mu^{1},\ldots,\mu^{m}),\,\nu=(\nu
^{1},\ldots,\nu^{m})\in (\fl^{*})^{m}$ we define semi-simple
orthogonal representations
\begin{itemize}
\item[]  $\rho^+_{\lambda}$ of $(\fl,\thl)$ on $\fa_+^{0,m}\oplus\fa_-^{0,m}
=\bigoplus_{i=1}^m \fa_+^{0,1}\oplus\fa_-^{0,1}$,
\item[]  $\rho^-_{\lambda}$ of $(\fl,\thl)$ on
$\fa_+^{m,0}\oplus\fa_-^{m,0}=\bigoplus_{i=1}^m \fa_+^{1,0}\oplus\fa_-^{1,0}$,
\item[]  $\tilde\rho^+_{\lambda}$ of $(\fl,\thl)$ on
$\fa_+^{m,0}\oplus\fa_-^{0,m}=\bigoplus_{i=1}^m \fa_+^{1,0}\oplus\fa_-^{0,1}$,
\item[]  $\tilde\rho^-_{\lambda}$ of $(\fl,\thl)$ on
$\fa_+^{0,m}\oplus\fa_-^{m,0}=\bigoplus_{i=1}^m \fa_+^{0,1}\oplus\fa_-^{1,0}$, and
\item[]  $\rho''_{\mu,\nu}$ of $(\fl,\thl)$ on
$\fa_+^{m,m}\oplus\fa_-^{m,m}=\bigoplus_{i=1}^m \fa_+^{1,1}\oplus\fa_-^{1,1}$
\end{itemize}
by
$$\rho^\pm_\lambda=\bigoplus_{i=1}^m \rho^\pm_{\lambda^i}\,,\quad
\tilde\rho^\pm_\lambda=\bigoplus_{i=1}^m \tilde\rho^\pm_{\lambda^i}\,,\quad
\rho''_{\mu,\nu}=\bigoplus_{i=1}^m \rho''_{\mu^i,\nu^i}.$$
Finally, let $(\rho_0)^{p,q}_{r,s}$ be the trivial representation of
$(\fl,\thl)$ on $\fa_+^{p,q}\oplus \fa_-^{r,s}$.

The symmetric group $\frak S_{m}$ acts on $(\fl^{*})^{m}$ by permuting
coordinates and on
$(\fl^{*})^{m}\oplus(\fl^{*})^{m}$ by permuting pairs of coordinates. The
group $(\ZZ_{2})^{m}$ acts on $(\fl^{*})^{m}$ by changing the signs of the
coordinates. We define the orbit spaces
$\Lambda_m:=(\fl^*\setminus  0)^{m}/\frak S_{m}\ltimes (\ZZ_{2})^{m}$ and
$\Lambda''_m:=\Big( ((\fl^*\setminus
0)^{m}/(\ZZ_{2})^{m})\oplus((\fl^*\setminus 0)^{m}/(\ZZ_{2})^{m}) \Big)/\frak
S_{m}$. Finally we define an action of $\Aut(\fl,\thl)$ on $\Lambda_m$ and $\Lambda''_m$
by
$S^*[\lambda]:=[S^*\lambda]$ and $S^*[\mu,\nu]:=[S^*\mu,S^*\nu]$.

\begin{pr} \label{PHom} Let $(\fl,\thl)$ be an abelian Lie algebra with
involution $\thl=-\Id$ and let $\fa=\fa_+^{p,q}\oplus\fa_-^{r,s}$ be as defined
above. Then
the map
 \begin{eqnarray*}
     \bigcup_{i \in I^{p,q}_{r,s}}
      \Lambda_{m_1}\times\dots\times
\Lambda_{m_4}\times\Lambda''_{m_5}&\longrightarrow&
      \Hom((\fl,\thl),\so(\fa))_{\rm ss} /\Aut(\fa)\\
     ([\lambda_1],\dots,[\lambda_4],[\mu,\nu])&\longmapsto&
      [\rho^+_{\lambda_1}\oplus \rho^-_{\lambda_2} \oplus
\tilde\rho^+_{\lambda_3} \oplus \tilde \rho^-_{\lambda_4}\
      \oplus    \rho''_{\mu,\nu}\oplus(\rho_0)^{p_0,q_0}_{r_0,s_0}]\,,
\end{eqnarray*}
where
$$I^{p,q}_{r,s}= \left\{i=(m_1,...,m_5,p_0,q_0,r_0,s_0)\in \ZZ^9
\left|\begin{array}{l} m_2+m_3+m_5+p_0=p,\\
m_1+m_4+m_5+q_0=q,\\m_2+m_4+m_5+r_0=r,\\ m_1+m_3+m_5+s_0=s
\end{array}\right.\right\}$$
is a
bijection. It is equivariant with respect to the action of $\Aut(\fl,\thl)$.
\end{pr}

\subsection{Admissible 3-dimensional Lie algebras with involution}

We know that any indecomposable, non-irreducible symmetric triple of
index~2 has the  structure of a quadratic extension $\dd$ in a canonical 
way, where  $(\fl,\thl)$ is an admissible Lie algebra with involution 
such that 
$\fl_{-}$  is at most 2-dimensional, $\fa$ is a suitable orthogonal 
$(\fl,\thl)$-module  and $[\alpha,\gamma]\in\cH^2_Q(\fl,\thl,\fa)_{0}$.
Therefore we  will now determine all admissible Lie algebras with involution 
$(\fl,\thl)$ such that 
$\fl_{-}$  is at most 2-dimensional. Moreover we 
determine   $\\cH^2_Q(\fl,\thl,\fa)_{0}$ for an arbitrary
semi-simple  orthogonal $(\fl,\thl)$-module~$\fa$.

\begin{pr}\label{P71}
If $(\fl,\theta_\fl)$ is an admissible Lie algebra with involution and
if ${\dim \fl_{-}\le 2}$, then $(\fl,\theta_\fl)$ is isomorphic to one of
the following Lie algebras with involution (given by the induced
decomposition into eigenspaces):
\begin{enumerate}
\item $\fl=\RR^k\ (k\le 2) ,\ \fl_+=0 ,\ \fl_-=\fl $,
\item $\fl=\fn(2)=\{[X,Y]=Z,\ [X,Z]=-Y\} ,\ \fl_+=\RR\cdot Z ,\
  \fl_-=\Span\{X,Y \}, $
\item $\fl=\fr_{3,-1}=\{[X,Y]=Y,\ [X,Z]=-Z\} ,\ \fl_+=\RR\cdot(Y+Z) ,\
  \fl_-=\Span\{X,\ Y-Z \}, $
\item $\fl=\fh(1)=\{[X,Y]=Z\} ,\ \fl_+=\RR\cdot Z ,\ \fl_-=\Span\{X,Y \} $
\item $\fl=\fsu(2)=\{[H,X]=2Y,\ [H,Y]=-2X,\ [X,Y]=2H\}\\ \fl_+= \RR\cdot H,\
\fl_-=\Span\{X,\ Y\},$
\item $\fl=\fsl(2,\RR)=\{[H,X]=2X,\ [H,Y]=-2Y,\ [X,Y]=H\}\\ \fl_+= \RR\cdot H,\
\fl_-=\Span\{X,\ Y\},$
\item $\fl=\fsl(2,\RR), \ \fl_+=\RR\cdot (X-Y),\ \fl_-=\Span\{H,\ X+Y \}.$
\end{enumerate}
\end{pr}
\proof
Let $(\fl,\theta_\fl)$ be a Lie algebra with involution satisfying 
$(T_{1})$.
Obviously, if $\dim \fl_{-}=1$, then $\fl=\fl_{-}=\RR$, and if $\dim 
\fl_{-}=2$ and $[\fl_{-},\fl_{-}]=0$, then $\fl=\fl_{-}=\RR^{2}$. 

If $\dim \fl_{-}=2$ and $[\fl_{-},\fl_{-}]\not=0$, then $(\fl,\theta_\fl)$ is 
isomorphic to one of the Lie algebras with involution 
$\fl_{\ph}$, $\ph\in \fsl(2,\RR)$, which are defined as follows:  
\begin{eqnarray*}
    &\fl_{\ph}=(\fl_{\ph})_{+}\oplus (\fl_{\ph})_{-},\ 
(\fl_{\ph})_{-}= \Span\{L_{1},L_{2}\}=\RR^{2},\ 
(\fl_{\ph})_{+}=\RR\cdot L_{3}=\RR&\\
&[L_{1},L_{2}]=L_{3},\ [L_{3},L]=\ph(L)\mbox{ for all } L\in 
(\fl_{\ph})_{-}.&
\end{eqnarray*}

For $\ph,\ph'\in \fsl(2,\RR)$ the Lie algebras with involution 
$\fl_{\ph}$ and $\fl_{\ph'}$ are isomorphic if and only if there is a map 
$g\in \GL(2,\RR)$ such that $\ph'=\det(g^{-1}) g\ph g^{-1}$. 
Consequently, $(\fl,\theta_\fl)$ is isomorphic to 
$\fl_{i}:=\fl_{\ph_{i}}$ for 
exactly one of the following $\ph_{i}\in\fsl(2,\RR)$, $i=1,\ldots,6$:
$$\ph_{1}=0,\ \ph_{2,3}=\small{\pm\left(
\begin{array}{cc}
    0 & 1  \\
    0 & 0
\end{array}
\right)},\ 
\ph_{4}=\small{\left(
\begin{array}{cc}
    1 & 0  \\
    0 & -1
\end{array}
\right)},\ 
\ph_{5,6}=\small{\pm\left(
\begin{array}{cc}
    0 & -1  \\
    1 & 0
\end{array}
\right)}.
$$
It is not hard to see that $\fl_{1}\cong\fh(1)$, 
$\fl_{2}\cong\fr_{3,-1}$, 
$\fl_{3}\cong\fn(2)$, 
$\fl_{5}\cong\fsu(2)$, all with the involution given in the 
proposition, $\fl_{4}\cong \fsl(2,\RR)$ with the involution given in {\it 
6.} and   $\fl_{6}\cong\fsl(2,\RR)$ with the involution given in {\it 
7.} We will see in this subsection that all these Lie algebras with 
involution are indeed admissible.
\qed

Next we will determine $\cH^2_Q(\fl,\thl,\fa)_{0}$ for all $(\fl,\thl)$ listed in
Prop.~\ref{P71} and any semi-simple orthogonal $(\fl,\thl)$-module~$\fa$. Let $\Theta$
be defined as in
(\ref{E9}). By Prop.~\ref{P31} we can identify $\cH^2_Q(\fl,\thl,\fa)$ with 
$\cH^2_Q(\fl,\fa)^\Theta$.
Therefore let us recall the following fact on quadratic cohomology sets of
three-dimensional Lie algebras (see \cite{KO2}, Lemma 7.2):

If $\fl$ is a three-dimensional  unimodular Lie algebra and $\fa$ is
an orthogonal $\fl$-module, then
\begin{eqnarray*}
\iota_Q: {\cal H}^{2}_Q(\fl,\fa) &\longrightarrow &
(H^2(\fl,\fa)\setminus\{ 0\})\cup C^3(\fl)\\
\,[\alpha,\gamma] &\longmapsto &\left\{
\begin{array}{ll} [\alpha]\in
      H^2(\fl,\fa)&\mbox{if } [\alpha]\not=0\\
\gamma\in C^3(\fl) & \mbox{if } [\alpha]=0
\end{array} \right.
\end{eqnarray*}
is a bijective map.

Now let $\theta_\fl$ be an involution on $\fl$ and let $\fa$ be
an orthogonal $(\fl,\thl)$-module.
 Obviously, $\Theta\circ \iota_Q=\iota_Q \circ \Theta$. In
particular, restricting  $\iota_Q$ to  ${\cal H}^{2}_Q(\fl,\fa)^\Theta$
we obtain a bijection
\begin{equation}\label{Eioq}
    \iota_Q~:\ {\cal H}^{2}_Q(\fl,\fa)^\Theta \longrightarrow  (H^2(\fl,\fa)^\Theta
\setminus\{ 0\})\cup C^3(\fl)_+.
\end{equation}
Here $H^2(\fl,\fa)^\Theta$ denotes the space of invariants of the action of $\Theta$ on 
$H^2(\fl,\fa)$.

\subsubsection{The case $\fl=\fn(2)$ or $\fl=\fr_{3,-1}$}

Now let $(\fl,\thl)$ be a three-dimensional admissible Lie algebra with
involution such that $\dim\fl'=2$, i.e. $\fl=\fn(2)$ or $\fl=\fr_{3,-1}$ with
$\thl$ as given in Proposition \ref{P71}.
The adjoint representation of $\fl$ induces a semi-simple
representation $\ad_0$ of $\fl_0:=\fl/\fl'$ on $\fl'$.
Let $\lambda ^i\in(\fl_{\Bbb C}/\fl'_{\Bbb C})^*$, $i=1,2$, be the weights of the complexification of $\ad_0$, and let
$$V_{\lambda^i}=\{U\in \fl'_{\Bbb C}\mid \ad_0(L)(U)=\lambda^i(L)\cdot U
\mbox{ for all  } L\in\fl\} $$
be the corresponding weight spaces.
Let $(\rho,\fa)$ be a semi-simple orthogonal
$(\fl,\theta_\fl)$-module. Then $\fl'\subset \ker\rho$ since $\fl'=R(\fl)$. Hence,
$\rho$ induces a representation of $\fl/\fl'$ on $\fa$. In particular,
the complexified module $(\rho,\fa_{\Bbb C})$ decomposes into
weight spaces
$$ E_\lambda=\{A\in\fa_{\Bbb C} \mid \rho(L)(A)=\lambda(L)\cdot A
\mbox{ for all } L\in\fl\},$$
where $\lambda\in (\fl_{\Bbb C}/\fl'_{\Bbb C})^*$. We will identify $\lambda$ with 
$\lambda(X)\in\CC$. For $\fl=\fn(2)$ we also define $E:=E_i\oplus E_{-i}$.

\begin{pr} \label{Pn2}
Suppose $\fl\in\{\fr_{3,-1},\ \fn(2)\}$. Let
$Z_{\fl}\subset C^2(\fl,\fa)_{+}$ be defined by
$$Z_{\fl}=\{\alpha\in C^2(\fl,\fa)\mid
\alpha(Y,Z)\in\fa ^\fl_+,\ \alpha(X,Y)\in E_1,\,
\tha(\alpha(X,Y))=-\alpha(X,Z)\}  $$
if $\fl=\fr_{3,-1}$ and by
$$Z_{\fl}=\{\alpha\in C^2(\fl,\fa)\mid
\alpha(Y,Z)\in\fa ^\fl_+,\ \alpha(X,Y)\in E\cap\fa_+,\,
\alpha(X,Z)=\rho(X)\alpha(X,Y)\}  $$
if $\fl=\fn(2)$. Then
\begin{eqnarray*}
(Z_{\fl}\setminus\{0\})\cup C^3(\fl)&\longrightarrow& {\cal
H}^{2}_Q(\fl,\thl,\fa) \\
(Z_{\fl}\setminus\{0\}) \ni \alpha & \longmapsto & [\alpha,0]\\
C^3(\fl)\ni\gamma&\longmapsto& [0,\gamma]
\end{eqnarray*}
is a bijection and
\begin{eqnarray*}
\, [\alpha,0]\in {\cal H}^{2}_Q(\fl,\thl,\fa)_0 &\Leftrightarrow&
\fa^\fl=\RR\cdot\alpha(Y,Z) \mbox{ and } \\&&\Span\{\alpha(X,Y), \alpha(X,Z)\}
\mbox{ is non-degenerate }\\[1ex]
\, [0,\gamma]\in {\cal H}^{2}_Q(\fl,\thl,\fa)_0&\Leftrightarrow& \fa^\fl=0
\mbox{ and } \gamma\not=0
\end{eqnarray*}
for all $\alpha\in Z_{\fl}\setminus\{0\}$ and $\gamma \in C^3(\fl)$.
\end{pr}
\proof
In \cite{KO2}, Lemma 7.1 we proved that
\begin{eqnarray}
&\{\alpha\in C^2(\fl,\fa)\mid \alpha(\fl',\fl')\subset
\fa^\fl,\ \alpha(X,V_{\lambda ^i})\subset
E_{\lambda ^i}, i=1,2\}  \longrightarrow  H^2(\fl,\fa)& \label{Eio}\\
&\alpha \longmapsto  [\alpha] &\nonumber
\end{eqnarray}
is correctly defined and an isomorphism.

By Equation (\ref{E10}) we have
$\rho(-X)(\theta_\fa(A))=\theta_\fa(\rho(X)(A))$ for all $A\in\fa$,
which implies
$$\tha (E_\lambda) = E_{-\lambda}.$$
Now let $\alpha\in C^2(\fl,\fa)$ be such that $ \alpha(X,V_{\lambda ^i})\subset
E_{\lambda ^i}$, $ i=1,2$. Then
$$(\Theta\alpha)(X,V_{\lambda ^i})=\tha (\alpha(X,\thl(V_{\lambda ^i})))=
  \tha (\alpha(X,V_{-\lambda ^i}))\subset \tha(E_{-\lambda
  ^i})=E_{\lambda ^i}.$$
Moreover, if $\alpha(\fl',\fl')\subset \fa ^\fl$, then
$$(\Theta\alpha)(\fl',\fl')=\tha (\alpha(\thl(\fl'),\thl(\fl'))
=\tha(\alpha(\fl',\fl'))\subset\fa ^\fl.$$
It follows that $\{\alpha\in C^2(\fl,\fa)\mid \alpha(\fl',\fl')\subset
\fa ^\fl,\ \alpha(X,V_{\lambda ^i})\subset
E_{\lambda ^i}, i=1,2\}  $ is $\Theta$-invariant. Hence the
restriction of the map defined in (\ref{Eio}) to
$$
   \{\alpha\in C^2(\fl,\fa)_+\mid  
\alpha(\fl',\fl')\subset
\fa ^\fl,\ \alpha(X,V_{\lambda ^i})\subset
E_{\lambda ^i}, i=1,2\} =Z_{\fl} $$ is an isomorphism onto
$H^2(\fl,\fa)^\Theta$. The first assertion of the proposition now
follows from (\ref{Eioq}) and $C^3(\fl)_+=C^3(\fl)$.
It remains to prove the statement on admissibility.

Assume first that
$[\alpha,0]\in {\cal H}^{2}_Q(\fl,\thl,\fa)_0$ for $\alpha\in
Z_{\fl}\setminus\{0\}$. Since $[\alpha,0]$ is indecomposable
we have $\fa^\fl=\RR\cdot\alpha(Y,Z)$. In particular,
$\RR\cdot\alpha(Y,Z)\subset \fa$ is a non-degenerate subspace. Assume that
$\Span\{\alpha(X,Y), \alpha(X,Z)\}\not=0$ and that $\ip_\fa$ restricted to this
subspace degenerates. Then $\Span\{\alpha(X,Y), \alpha(X,Z)\}\not=0$ is totally
isotropic. Let $\fb_{1}$ be as in Condition $(B_{1})$ of Definition
\ref{zwiesel}. Since obviously the submodule $\alpha(Y,Z)^\perp \cap \{\alpha(X,Y),
\alpha(X,Z)\}^\perp$ is contained in $\fb_1$ and $\fb_1$ is non-degenerate we
obtain $\alpha(Y,Z)^\perp \subset \fb_1$. This yields
$\alpha(X,Y)=\alpha(X,Z)=0$, a contradiction.

Now assume that $\alpha\in Z_{\fl}\setminus\{0\}$ satisfies
$\fa^\fl=\RR\cdot\alpha(Y,Z)$ and that the subspace $\Span\{\alpha(X,Y),
\alpha(X,Z)\}\subset \fa$ is non-degenerate. Since
$\alpha_0(\Ker\lb_\fl)=\RR\cdot\alpha(Y,Z)$ Condition $(B_0)$ is satisfied. We
prove that $(B_1)$ is also satisfied. We use that the submodule
$\fm:=\alpha(Y,Z)^\perp \cap \{\alpha(X,Y), \alpha(X,Z)\}^\perp =
(\fa^\fl)^\perp \cap \{\alpha(X,Y), \alpha(X,Z)\}^\perp$ is contained in $\fb_1$
and that $\fb_1\subset (\fa^\fl)^\perp$ itself is a submodule. In case
$\fl=\fn(2)$ this implies that $\fb_1$ is equal to $(\fa^\fl)^\perp$ or to
$\fm$. Since both submodules are non-degenerate $(B_1)$ holds. Now we consider
the case $\fl=\fr_{3,-1}$. Assume that $\fb_1$ contains the submodule $\RR\cdot
\alpha(X,Y)$. Then
$$\langle\alpha(X,Z),\alpha(X,Y)\rangle = \langle (\rho(X)+\Id)(\Phi(Z)),
\alpha(X,Y)\rangle =0$$
since $\alpha(X,Y)\in E_1$. Since on the other hand $\alpha(X,Z)\in E_{-1}$ and
by assumption  $\Span\{\alpha(X,Y), \alpha(X,Z)\}\subset \fa$ is non-degenerate
we obtain $\alpha(X,Y)=\alpha(X,Z)=0$. Similarly  $\RR\cdot \alpha(X,Z)\subset
\fb_1$ implies $\alpha(X,Y)=\alpha(X,Z)=0$. Again we obtain that $\fb_1$ is
equal to $(\fa^\fl)^\perp$ or to $\fm$ and $(B_1)$ holds.  Now let again
$\fl$ be $\fr_{3,-1}$ or $\fn(2)$. Since $R_2(\fl)=0$ Conditions $(A_k)$ and
$(B_k)$ hold for $k\ge 2$. Moreover, $(A_0)$ holds since $\fz(\fl)=0$ and
$(A_1)$ holds since no $\alpha\in Z_{\fl}\setminus\{0\}$
satisfies Assumption $(i)$ of $(A_1)$. Since
also Conditions $(T_1)$ and $(T_2)$ are satisfied $[\alpha,0]$ is admissible by
Theorem \ref{twain}. Since $\fl$ does not decompose into the direct sum of two
non-vanishing Lie algebras and $\fa^\fl=\RR\cdot\alpha(Y,Z)$ Proposition
\ref{hermlin} yields $[\alpha,0]\in {\cal H}^{2}_Q(\fl,\thl,\fa)_0$.

For $[0,\gamma]$ with  $\gamma\in C^3(\fl)$ Assumption (ii) of $(A_1)$ is
satisfied if and only if $\gamma=0$. Hence $(A_1)$ holds if and only if
$\gamma\not=0$. Since all other Conditions $(A_k)$ and $(B_k)$ are trivially
satisfied $[0,\gamma]$ is admissible if and only if $\gamma\not=0$ and
$\fa_+^\fl=0$ (Condition $(T_{2})$). Applying
Proposition \ref{hermlin} we obtain $[0,\gamma]\in {\cal
H}^{2}_Q(\fl,\thl,\fa)_0$ if and only if  $\fa^\fl=0$ and $\gamma\not=0$.
\qed

\begin{pr} \label{PAut}
If $\fl=\fn(2)$, then
$$\Aut(\fl,\thl)=\Big\{
{\small
\left( \begin{array}{ccc}
u&0&0\\
v&a&0\\
0&0&ua
\end{array}\right)
}
\ \Big|\ u=\pm1,\, a,v\in\RR,\ a\not=0\Big\}$$

For $\fl=\fr_{3,-1}$ we have
$$\Aut(\fl,\thl)=\Big\{
{\small
    \left( \begin{array}{ccc}
    1&0&0\\
    v&a&0\\
    -v&0&a
    \end{array}\right)}\ \Big|\ a,v\in \RR, a\not=0 \Big\} \cup \Big\{
    \small{
    \left( \begin{array}{ccc}
    -1&0&0\\
    v&0&b\\
    -v&b&0
    \end{array}\right)}\ \Big|\ b,v\in \RR, b\not=0\Big\}.$$
Here all automorphisms are written with respect to the basis $X,Y,Z$ of $\fl$.  
\end{pr}

\subsubsection{The case $\fl=\fh(1)$ }
Now we suppose that $\fl=\fh(1)$ and that $\thl$ is given as in Proposition
\ref{P71}, 4.

\begin{pr}
\label{Ph1}
For $\fl=\fh(1)$ and
$$Z_\fl:= \{\alpha\in C^2(\fl,\fa)\mid
\alpha(X,Y)=0
,\ \alpha(Z,\fl)\subset \fa^\fl_-\}  $$
the map
\begin{eqnarray*}
(Z_\fl\setminus\{0\})\cup C^3(\fl)&\longrightarrow& {\cal
H}^{2}_Q(\fl,\thl,\fa) \\
Z_\fl\setminus\{0\} \ni \alpha & \longmapsto & [\alpha,0]\\
C^3(\fl)\ni\gamma&\longmapsto& [0,\gamma]
\end{eqnarray*}
is a bijection and ${\cal H}^{2}_Q(\fl,\thl,\fa)_0$ is the image of
$$\{\alpha \in Z_\fl\setminus\{0\} \mid
\alpha(Z,\fl)=\fa^\fl\}$$
under this bijection. In particular, ${\cal H}^{2}_Q(\fl,\thl,\fa)_0
=\emptyset$ if $\fa^\fl_+\not=0$.
\end{pr}
\proof
In \cite{KO2}, Lemma 7.5 we proved that
\begin{eqnarray*}
&\{\alpha\in C^2(\fl,\fa)\mid \alpha(X,Y)=0
,\ \alpha(Z,\fl)\subset \fa^\fl\}  \longrightarrow  H^2(\fl,\fa)&\\
&\alpha \longmapsto  [\alpha] &
\end{eqnarray*}
is an isomorphism. The domain of this isomorphism is $\Theta$-invariant. Hence
it restricts to an isomorphism from $\{\alpha\in C^2_+(\fl,\fa)\mid
\alpha(X,Y)=0
,\ \alpha(Z,\fl)\subset \fa^\fl\}=Z_\fl$ to $H^2(\fl,\fa)^\Theta$.

The cohomology classes $[0,\gamma]\in {\cal H}^{2}_Q(\fl,\thl,\fa)$ are
not admissible, neither $(A_0)$ nor $(A_1)$ is satisfied (see also \cite{KO2}, Lemma 
7.6). If
$\alpha \in Z_\fl\setminus\{0\}$, then neither assumption
$(i)$ of $(A_0)$ nor assumption $(i)$ of $(A_1)$ is satisfied. Hence $(A_0)$ and
$(A_1)$ hold. Furthermore, both $(B_0)$ and $(B_1)$ are equivalent to the
condition that $\alpha(Z,\fl)$ is non-degenerate. Condition $(T_{2})$ is equivalent
to $\fa^\fl_+=0$. Hence, $[\alpha,0]$ is admissible if and only if $\fa^\fl_+=0$
and $\alpha(Z,\fl)$ is non-degenerate. Proposition \ref{hermlin} now yields
$[\alpha,0]\in {\cal H}^{2}_Q(\fl,\thl,\fa)_0$ if and only if
$\alpha(Z,\fl)=\fa^\fl$.
\qed
\begin{pr}\label{PAut2}
For $\fl=\fh(1)$ we have
$$\Aut(\fl,\thl)=\Big\{
    \small{
\left( \begin{array}{cc}
A&0\\
0&u
\end{array}\right)
}
\ \Big|\ A\in GL(2,\RR),\ \det A=u\
\Big\}, $$
where the automorphisms are written with respect to the basis $X,Y,Z$ of $\fl$.
\end{pr}

\subsubsection{The case $\fl=\RR^k$, $k=1,2$}  
\begin{pr}     If $\fl=\RR^k$, $k=1,2$, then we can identify
$$ \cH^2_Q(\fl,\thl,\fa)=H^2(\fl,\fa)^{\Theta}=C^2_{+}(\fl,\fa^\fl)=
C^2(\fl,\fa^\fl_{+}).$$ and we have    
$$\cH^2_Q(\fl,\thl,\fa)_0=    
\left\{
        \begin{array}{ll}            
        C^2(\fl,\fa^\fl_{+})\setminus \{0\}&\mbox{ if }\ \dim            
        \fa^\fl_{+}=1,\ \dim \fa^\fl_{-}=0\\            
        \{0\} &\mbox{ if }\ \dim \fa^\fl=0, \ (\fl,\thl,\fa)\mbox{ indecomposable}\\              
        \emptyset &\mbox{ otherwise .}         \end{array}     \right.$$
\end{pr}
Here the triple $(\fl,\thl,\fa)$ is indecomposable, if it has not any non-trivial 
decomposition in the sense of Definition \ref{kosel}. The proof of this proposition is 
easy, so we will omit it.

\subsubsection{The case $\fl=\fsu(2)$ or $\fsl(2,\RR)$}

\begin{lm}\label{whitehead} Let $\fl\in\{\fsu(2),\ \fsl(2,\RR)\}$. Let $\theta_\fl$ be 
an involution of $\fl$.
Then we have for all semi-simple orthogonal $(\fl,\theta_\fl)$-modules $\fa$
$$ \cH^2_Q(\fl,\thl,\fa)= C^3(\fl)\ .$$
Moreover,
$$\cH^2_Q(\fl,\thl,\fa)_0=    
\left\{
        \begin{array}{ll}            
        C^3(\fl)&\mbox{ if }\ \fa^\fl=0\\            
        \emptyset &\mbox{ if }\ \fa^\fl\ne 0
\end{array}     \right.\ .$$ 
\end{lm}

\proof
Since $\fl$ is semi-simple we have $H^2(\fl,\fa)=0$. In order to obtain the first 
assertion we now
combine Proposition \ref{P31} with (\ref{Eioq}). Observe that $C^3(\fl)_+=C^3(\fl)$. The 
second assertion
is then easy to check.
\qed

Next we introduce certain orthogonal $(\fl,\theta_\fl)$-modules for
$\fl\in\{\fsu(2),\fsl(2,\RR)\}$. Since we are interested in admissible quadratic 
extensions of $(\fl,\thl)$ 
by such orthogonal 
$(\fl,\theta_\fl)$-modules $\fa$ which yield symmetric triples of index 2 and since here 
$\dim \fl_-=2$ 
we restrict ourselves to those $\fa$ for which $\ip_\fa$ restricted to $\fa_-$ is 
positive definite 
(see Remark~\ref{R51}).

Let $\fl=\fsu(2)$. For $k\in\NN$ let $\rho^\pm_k$
be the irreducible $\fl$-representation on a real vector space $\fa$ of
dimension $2k+1$ which is equipped
with a positive definite $\fl$-invariant scalar product and an involution $\tha$
uniquely characterized by (\ref{E10}) and
$\tha|_{\fa^{\fl_+}}=\pm(-1)^k\Id_{\fa^{\fl_+}}$. Then $\rho_k^+$
acts on $\fa_+^{0,k+1}\oplus\fa_-^{0,k}$, and $\rho_k^-$
acts on $\fa_+^{0,k}\oplus\fa_-^{0,k+1}$. By $\rho'_k$ we denote
the irreducible orthogonal representation of $\fl$ on a real vector space $\fa$
of dimension $4k$. Then $\fa$ carries an isometric involution $\tha$ satisfying
$(\ref{E10})$ which is uniquely determined up to a sign. A different choice
of the sign would produce an equivalent orthogonal $(\fl,\theta_\fl)$-module.
Thus there is no need to fix it here. We have $\fa=\fa_+^{0,2k}\oplus\fa_-^{0,2k}$.

Let $\fl=\fsl(2,\RR)$, and let $\thl$ as in Proposition \ref{P71},7.
For $k\in\NN$ let $\rho^\pm_k$
be the irreducible $\fl$-representation on a real vector space $\fa$ of
dimension $2k+1$ which is equipped
with an involution $\tha$
uniquely characterized by (\ref{E10}) and
$\tha|_{\fa^{\fl_+}}=\pm(-1)^k\Id_{\fa^{\fl_+}}$.
Then there is an $\fl$-invariant scalar product on $\fa$ which is positive
definite on $\fa_-$. Then $\rho_k^+$
acts on $\fa_+^{k+1,0}\oplus\fa_-^{0,k}$, and $\rho_k^-$
acts on $\fa_+^{k,0}\oplus\fa_-^{0,k+1}$.
We also consider the real irreducible  $\fl$-representation acting on a
$2k$-dimensional real vector space $V_k$. The
natural $\fl$-representation $\rho'_k$ on $\fa:=V_k\oplus V_k^*$ carries an
invariant scalar product of signature $(2k,2k)$ induced by the dual pairing.
Moreover, there is exactly one involution $\theta_\fa$ which satisfies
(\ref{E10}) and has a positive definite $(-1)$-eigenspace. Note that $\tha$ switches
the two summands $V_k$ and $V_k^*$.
We have $\fa=\fa_+^{2k,0}\oplus \fa_-^{0,2k}$.

For $\fl=\fsl(2,\RR)$ and $\thl$ as in Proposition \ref{P71},6.\ we only
consider the orthogonal $(\fl,\thl)$-module given by
$(\rho_1^+,\fa,\ipa,\tha):=(\ad,\fl,B, -\thl)$, where $B$ denotes the
Killing form on $\fl$. We have
$\fa=\fa_+^{1,1}\oplus \fa_-^{0,1}$.

We call an orthogonal $(\fl,\thl)$-module $\fa$ indecomposable if it has no
proper non-trivial non-degenerate $(\fl,\thl)$-invariant submodule.
A semi-simple indecomposable module
is either irreducible or the direct
sum of two irreducible totally isotropic $(\fl,\thl)$-modules in duality.
The latter case does not occur if the restriction of the scalar product to
$\fa_+$ or $\fa_-$ is definite. All the $(\fl,\thl)$-modules defined above are
irreducible.

\begin{lm}\label{rep}
Let $(\fl,\thl)$ be as in cases 5.-7. of Proposition \ref{P71}.
Then the orthogonal $(\fl,\thl)$-modules just defined
exhaust the equivalence classes of those indecomposable orthogonal
$(\fl,\thl)$-modules satisfying $\fa^\fl=\{0\}$ and $\fa_-$
is positive definite.
They
are pairwise inequivalent.
\end{lm}
\proof
See \cite{CP}, Ch.V,\S 3.
\qed

If $k=(k_1,\dots,k_p)\in \NN^p$ for some $p\ge 0$, then we denote by
$\rho_k^+$ the direct sum module $\rho^+_{k_1}\oplus\dots\oplus \rho^+_{k_p}$.
By convention $\NN^0=\emptyset$, and the corresponding direct
sum is the zero module. In the same way we define $\rho_k^-$ and $\rho'_k$,
$k\in\NN^p$. Moreover we set $|k|:=k_1+\dots+k_p$.

\subsection{The classification result}

Now we can formulate our classification of symmetric triples of index 2. As above we 
will use the notation
$\fl_0$ for $\fl/R(\fl)$ for a given solvable Lie algebra $\fl$. Furthermore let
$\bar{\frak S}_{p,q}$ be the group 
$(\frak S_{p}\ltimes(\ZZ_{2})^{p})\times(\frak{S}_{q}\ltimes(\ZZ_{2})^{q})$ 
which acts on 
$(\fl_0^*\setminus 0)^p\times(\fl_0^*\setminus 0)^q$.

\begin{theo}\label{Klasse}
If $(\fg,\ip,\theta)$ is a symmetric triple associated
with an indecomposable non-semi-simple symmetric
space of index 2, then it is isomorphic to
$\dd$ for exactly one of the data in the following list (which contains only data giving 
rise to such triples):
\begin{enumerate}
\item $\fl=\RR^1=\RR\cdot X,\ \fl_+=0 ,\ \fl_-=\fl $,
\begin{enumerate}
\item $\fa=\fa_+^{0,1}\oplus\fa_-^{1,0}  \oplus \fa_+^{p,q}\oplus\fa_-^{0,p+q}$,
$\ p,q\ge 0,$\\
$\rho=\tilde\rho_1^-  \oplus \tilde\rho_\lambda^+\oplus\rho^+_\mu$, \\
$ \lambda=(\lambda^1,\dots,\lambda^p)\in (\fl^*)^p=\RR^p,\
\mu=(\mu^1,\dots,\mu^q)\in (\fl^*)^q=\RR^q$, \\
$0<\lambda^1\le\dots \le \lambda^p$, $0<\mu^1\le\dots \le \mu^q$,\\
$\alpha=0,\ \gamma=0$;
\item $\fa=\fa_+^{1,0}\oplus\fa_-^{1,0}  \oplus \fa_+^{p,q}\oplus\fa_-^{0,p+q}$,
$\ p,q\ge 0,$\\
$\rho=\rho_1^-  \oplus \tilde\rho_\lambda^+\oplus\rho^+_\mu$\\
$ \lambda=(\lambda^1,\dots,\lambda^p)\in (\fl^*)^p=\RR^p,\
\mu=(\mu^1,\dots,\mu^q)\in (\fl^*)^q=\RR^q$, \\
$0<\lambda^1\le\dots \le \lambda^p$, $0<\mu^1\le\dots \le \mu^q$,\\$\alpha=0, \
\gamma=0$;
\item $\fa=\fa_+^{1,1}\oplus\fa_-^{1,1}  \oplus \fa_+^{p,q}\oplus\fa_-^{0,p+q}$,
$\ p,q\ge 0,$\\
$\rho=\rho''_{1,\nu}  \oplus \tilde\rho_\lambda^+\oplus\rho^+_\mu$\\
$\nu\in\fl^*\setminus 0=\RR^1\setminus 0,\ \lambda=(\lambda^1,\dots,\lambda^p)\in 
(\fl^*)^p=\RR^p,\
\mu=(\mu^1,\dots,\mu^q)\in (\fl^*)^q=\RR^q$, \\
$0<\lambda^1\le\dots \le \lambda^p$, $0<\mu^1\le\dots \le \mu^q$,
\\$\alpha=0,\ \gamma=0$;
\end{enumerate}
\item $\fl=\RR^2=\Span\{Y,Z\} ,\ \fl_+=0 ,\ \fl_-=\fl $,
\begin{enumerate}
\item $\fa= \fa_+^{p,q}\oplus\fa_-^{0,p+q}$, $\ p,q\ge 0,$ $p+q\ge3,$\\
$\rho=\tilde\rho_\lambda^+\oplus\rho^+_\mu$, \\
$ [\lambda,\mu]\in ((\fl^*\setminus 0)^p\times (\fl^*\setminus 0)^q)/_\sim\,,$
such that
\begin{center}$\{(\lambda^i(Y),\lambda^i(Z))\mid i=1,\dots,p\}\ \cup\
\{(\mu^j(Y),\mu^j(Z))\mid j=1,\dots,q\}$
\end{center}
is not contained in the union of two one-dimensional subspaces of $\RR^2$,\\
and\\
$(\lambda_1,\mu_1)\sim (\lambda_2,\mu_2)\ \Leftrightarrow\ $
$\Span\{y_1,z_1\}=\Span\{y_2,z_2\}\ \mod\ \bar{\frak S}_{p,q}$,\\
where $y_i:=(\lambda_i^1(Y),\dots,\lambda_i^p(Y),\mu_i^1(Y),\dots,\mu_i^q(Y))$ and\\
$z_i:=(\lambda_i^1(Z),\dots,\lambda_i^p(Z),\mu_i^1(Z),\dots,\mu_i^q(Z))$ for
$i=1,2$,\\
$\alpha=0,\ \gamma=0$;\\
\item $\fa=\fa_+^{0,1}  \oplus \fa_+^{p,q}\oplus\fa_-^{0,p+q}$, $\ p,q\ge 0,$\\
$\rho=(\rho_0)^{0,1}_{0,0}  \oplus \tilde\rho_\lambda^+\oplus\rho^+_\mu$, \\
$[ \lambda,\mu]\in ((\fl^*\setminus 0)^p\times  (\fl^*\setminus 0)^q)/_\sim,$
with\\
$(\lambda_1,\mu_1)\sim (\lambda_2,\mu_2)\ \Leftrightarrow $
\begin{center}
$ (\Span\{y_1,z_1\},y_1\wedge z_1)=(\Span\{y_2,z_2\},\pm y_2\wedge z_2)\ \mod \
\bar{\frak S}_{p,q}$,
\end{center}
where $y_i, z_i$, $i=1,2$, are as in 2.\,(a),\\
$\alpha(Y,Z)=A_0$, where $A_0$ is the standard basis of $\fa_+^{0,1}$,\\
$\gamma=0$; \\
\item $\fa=\fa_+^{1,0}  \oplus \fa_+^{p,q}\oplus\fa_-^{0,p+q}$, $\ p,q\ge 0,$\\
$\rho=(\rho_0)^{1,0}_{0,0}  \oplus \tilde\rho_\lambda^+\oplus\rho^+_\mu$, \\
$[ \lambda,\mu]\in ((\fl^*\setminus 0)^p\times (\fl^*\setminus 0)^q)/_\sim,$ with\\
$(\lambda_1,\mu_1)\sim (\lambda_2,\mu_2)\ \Leftrightarrow $
\begin{center}
$ (\Span\{y_1,z_1\},y_1\wedge z_1)=(\Span\{y_2,z_2\},\pm y_2\wedge z_2)\ \mod \
\bar{\frak S}_{p,q}$,
\end{center}
where $y_i, z_i$, $i=1,2$, are as in 2.\,(a),\\
 $\alpha(Y,Z)=A_0$, where $A_0$ is the standard basis of $\fa_+^{1,0}$,\\
$\gamma=0$;
\end{enumerate}
\item $\fl=\fn(2),\ \fl_+=\RR\cdot Z ,\ \fl_-=\Span\{X,Y \}, $
\begin{enumerate}
\item $\fa= \fa_+^{p,q}\oplus\fa_-^{0,p+q}$, $\ p,q\ge 0,$\\
$\rho=\tilde\rho_\lambda^+\oplus\rho^+_\mu$, \\
$ \lambda=(\lambda^1,\dots,\lambda^p)\in (\fl_0^*)^p=\RR^p,\
\mu=(\mu^1,\dots,\mu^q)\in (\fl_0^*)^q=\RR^q$, \\
$0<\lambda^1\le\dots \le \lambda^p$, $0<\mu^1\le\dots \le \mu^q$, \\
$\alpha=0,\ \gamma(X,Y,Z)=\kappa, \ \kappa\in\{1,-1\}$;
\item $\fa=\fa_-^{0,1}\oplus \fa_+^{p,q}\oplus\fa_-^{0,p+q}$, $\ p,q\ge 0,$\\
$\rho=(\rho_0)^{0,0}_{0,1}\oplus \tilde\rho_\lambda^+\oplus\rho^+_\mu$, \\
$ \lambda=(\lambda^1,\dots,\lambda^p)\in (\fl_0^*)^p=\RR^p,\
\mu=(\mu^1,\dots,\mu^q)\in (\fl_0^*)^q=\RR^q$, \\
$0<\lambda^1\le\dots \le \lambda^p$, $0<\mu^1\le\dots \le \mu^q$, \\
$\alpha(Y,Z)=A_0$, where $A_0$ is the standard basis of $\fa_-^{0,1}$,\\
$\gamma=0$;
\item $\fa=\fa_+^{0,1}\oplus\fa_-^{0,1}\oplus \fa_+^{p,q}\oplus\fa_-^{0,p+q}$,
$\ p,q\ge 0,$\\
$\rho=\rho^+_1\oplus \tilde\rho_\lambda^+\oplus\rho^+_\mu$, \\
$ \lambda=(\lambda^1,\dots,\lambda^p)\in (\fl_0^*)^p=\RR^p,\
\mu=(\mu^1,\dots,\mu^q)\in (\fl_0^*)^q=\RR^q$, \\
$0<\lambda^1\le\dots \le \lambda^p$, $0<\mu^1\le\dots \le \mu^q$, \\
$\alpha(X,Y)=A_1$, $\alpha(X,Z)=A_2$, $\alpha(Y,Z)=0$,\\ where $A_1,A_2$ is the
standard basis of $\fa_+^{0,1}\oplus\fa_-^{0,1}$,\\
$\gamma=0$;
\item $\fa=\fa_-^{0,1}\oplus\fa_+^{0,1}\oplus\fa_-^{0,1}\oplus
\fa_+^{p,q}\oplus\fa_-^{0,p+q}$, $\ p,q\ge 0,$\\
$\rho=(\rho_0)^{0,0}_{0,1}\oplus  \rho^+_1\oplus
\tilde\rho_\lambda^+\oplus\rho^+_\mu$, \\
$ \lambda=(\lambda^1,\dots,\lambda^p)\in (\fl_0^*)^p=\RR^p,\
\mu=(\mu^1,\dots,\mu^q)\in (\fl_0^*)^q=\RR^q$, \\
$0<\lambda^1\le\dots \le \lambda^p$, $0<\mu^1\le\dots \le \mu^q$, \\
$\alpha(X,Y)=rA_1$, $\alpha(X,Z)=rA_2$, $\alpha(Y,Z)=A_0$, $r\in\RR$, $r>0$,\\
where $A_0,A_1,A_2$ is the standard basis of
$\fa_-^{0,1}\oplus\fa_+^{0,1}\oplus\fa_-^{0,1}$,\\
$\gamma=0$;
\end{enumerate}
\item $\fl=\fr_{3,-1},\ \fl_+=\RR\cdot(Y+Z) ,\ \fl_-=\Span\{X,\ Y-Z \}, $
\begin{enumerate}
\item $\fa= \fa_+^{p,q}\oplus\fa_-^{0,p+q}$, $\ p,q\ge 0,$\\
$\rho=\tilde\rho_\lambda^+\oplus\rho^+_\mu$, \\
$ \lambda=(\lambda^1,\dots,\lambda^p)\in (\fl_0^*)^p=\RR^p,\
\mu=(\mu^1,\dots,\mu^q)\in (\fl_0^*)^q=\RR^q$, \\
$0<\lambda^1\le\dots \le \lambda^p$, $0<\mu^1\le\dots \le \mu^q$, \\
 $\alpha=0,\ \gamma(X,Y,Z)=\kappa, \ \kappa\in\{1,-1\}$;
\item $\fa=\fa_-^{0,1}\oplus \fa_+^{p,q}\oplus\fa_-^{0,p+q}$, $\ p,q\ge 0,$\\
$\rho=(\rho_0)^{0,0}_{0,1}\oplus \tilde\rho_\lambda^+\oplus\rho^+_\mu$, \\
$ \lambda=(\lambda^1,\dots,\lambda^p)\in (\fl_0^*)^p=\RR^p,\
\mu=(\mu^1,\dots,\mu^q)\in (\fl_0^*)^q=\RR^q$, \\
$0<\lambda^1\le\dots \le \lambda^p$, $0<\mu^1\le\dots \le \mu^q$, \\
$\alpha(Y,Z)=A_0$, where $A_0$ is the standard basis of $\fa_-^{0,1}$,\\
$\gamma=0$;
\item $\fa=\fa_+^{1,0}\oplus\fa_-^{0,1}\oplus \fa_+^{p,q}\oplus\fa_-^{0,p+q}$,
$\ p,q\ge 0,$\\
$\rho=\tilde\rho^+_1\oplus \tilde\rho_\lambda^+\oplus\rho^+_\mu$, \\
$ \lambda=(\lambda^1,\dots,\lambda^p)\in (\fl_0^*)^p=\RR^p,\
\mu=(\mu^1,\dots,\mu^q)\in (\fl_0^*)^q=\RR^q$, \\
$0<\lambda^1\le\dots \le \lambda^p$, $0<\mu^1\le\dots \le \mu^q$, \\
$\alpha(X,Y-Z)=A_1$, $\alpha(X,Y+Z)=A_2$, $\alpha(Y,Z)=0$,\\ where $A_1,A_2$ is
the standard basis of $\fa_+^{1,0}\oplus\fa_-^{0,1}$,\\
$\gamma=0$;
\item $\fa=\fa_-^{0,1}\oplus\fa_+^{1,0}\oplus\fa_-^{0,1}\oplus
\fa_+^{p,q}\oplus\fa_-^{0,p+q}$, $\ p,q\ge 0,$\\
$\rho=(\rho_0)^{0,0}_{0,1}\oplus  \tilde\rho^+_1\oplus
\tilde\rho_\lambda^+\oplus\rho^+_\mu$, \\
$ \lambda=(\lambda^1,\dots,\lambda^p)\in (\fl_0^*)^p=\RR^p,\
\mu=(\mu^1,\dots,\mu^q)\in (\fl_0^*)^q=\RR^q$, \\
$0<\lambda^1\le\dots \le \lambda^p$, $0<\mu^1\le\dots \le \mu^q$, \\
$\alpha(X,Y-Z)=rA_1$, $\alpha(X,Y+Z)=rA_2$, $\alpha(Y,Z)=A_0$, $r\in\RR$,
$r>0$,\\ where $A_0,A_1,A_2$ is the standard basis of
$\fa_-^{0,1}\oplus\fa_+^{1,0}\oplus\fa_-^{0,1}$,\\
$\gamma=0$;
\end{enumerate}
\item $\fl=\fh(1),\ \fl_+=\RR\cdot Z ,\ \fl_-=\Span\{X,Y \} $
\begin{enumerate}
\item $\fa=\fa_-^{0,1}\oplus \fa_+^{p,q}\oplus\fa_-^{0,p+q}$, $\ p,q\ge 0,$\\
$\rho=(\rho_0)^{0,0}_{0,1}\oplus \tilde\rho_\lambda^+\oplus\rho^+_\mu$, \\
$ [\lambda,\mu]\in ((\fl_0^*\setminus 0)^p\times\ (\fl_0^*\setminus
0)^q)/_\sim\,,$ with\\
$(\lambda_1,\mu_1)\sim (\lambda_2,\mu_2)\ \Leftrightarrow $ $(\lambda_1,\mu_1)$,
$(\lambda_2,\mu_2)$ satisfy $(i)$ or $(ii)$:
        \begin{itemize}
        \item[(i)]
        $\dim \Span \{x_1,y_1\}=\dim\Span\{x_2,y_2\}=1$ and
        $$(\Span \{x_1,y_1\},\RR\cdot y_1) =
          (\Span\{x_2,y_2\},\RR\cdot y_2)\mod\ \bar{\frak S}_{p,q},$$
        \item[(ii)]
        $\dim \Span \{x_1,y_1\}=\dim \Span\{x_2,y_2\}=2$
        and
        $$\exists\, r\in\RR,r\not=0\,:\ (x_2\wedge y_2,y_2)=(rx_1\wedge
        y_1,r^2y_1)\mod\ \bar{\frak S}_{p,q},$$
        \end{itemize}
where $x_i:=(\lambda_i^1(X),\dots,\lambda_i^p(X),\mu_i^1(X),\dots,\mu_i^q(X))$ and\\
$y_i:=(\lambda_i^1(Y),\dots,\lambda_i^p(Y),\mu_i^1(Y),\dots,\mu_i^q(Y))$ for
$i=1,2$,\\        
$\alpha(X,Z)=A_0$, $\alpha(X,Y)=\alpha(Y,Z)=0$,\\ where $A_0$ is the standard
basis of $\fa_-^{0,1}$,\\
$\gamma=0$;
\item $\fa=\fa_-^{0,2}\oplus \fa_+^{p,q}\oplus\fa_-^{0,p+q}$, $\ p,q\ge 0,$\\
$\rho=(\rho_0)^{0,0}_{0,2}\oplus \tilde\rho_\lambda^+\oplus\rho^+_\mu$, \\
$ [\lambda,\mu]\in ((\fl_0^*\setminus 0)^p\times\ (\fl_0^*\setminus
0)^q)/_\sim\,,$ with\\
$(\lambda_1,\mu_1)\sim (\lambda_2,\mu_2)\ \Leftrightarrow M_1 M_1^\top =
M_2M_2^\top \ \mod\ \bar{\frak S}_{p,q},$\\
for the $((p+q)\times 2)$-matrices $M_i:=(x_i^\top,y_i^\top)$, $i=1,2$, where $x_i,
y_i$ are defined as in 5.\,(a), \\
$\alpha(X,Z)=A_1$, $\alpha(Y,Z)=A_2$, $\alpha(X,Y)=0$,\\ where $A_1,A_2$ is the
standard basis of $\fa_-^{0,2}$,\\
$\gamma=0$;\end{enumerate}
\item $\fl=\fsu(2),\ \fl_+= \RR\cdot H,\ \fl_-=\Span\{X,\ Y\},$
\begin{enumerate}
\item[] $\fa=\fa_+^{0,|k|+|l|+2|m|+p}\oplus\fa_-^{0,|k|+|l|+2|m|+q}$, $\
p,q,r\ge 0$, \\
$k\in \NN^p,\ l\in\NN^q,\  m\in \NN^r$, \\
$k_1\le k_2\le\dots\le k_p,\ l_1\le\dots\le l_q,\ m_1\le\dots\le m_r$, \\
$\rho=\rho_k^+\oplus\rho_l^-\oplus\rho'_m$, \\
$\alpha=0$,\\
$\gamma(H,X,Y)=c,\ c\in\RR$.
\end{enumerate}
\item $\fl=\fsl(2,\RR),\ \fl_+= \RR\cdot H,\ \fl_-=\Span\{X,\ Y\},$
\begin{enumerate}
\item[] $\fa=\fa_+^{p,p}\oplus\fa_-^{0,p}$, $\ p\ge 0$,
$\ \displaystyle \rho=\bigoplus_{i=1}^p\rho_1^+$, \\
$\alpha=0$,
$\gamma(H,X,Y)=c,\ c\in\RR$.
\end{enumerate}
\item $\fl=\fsl(2,\RR), \ \fl_+=\RR\cdot (X-Y),\ \fl_-=\Span\{H,\ X+Y \}.$
\begin{enumerate}
\item[] $\fa=\fa_+^{|k|+|l|+2|m|+p,0}\oplus\fa_-^{0,|k|+|l|+2|m|+q}$, $\
p,q,r\ge 0$, \\
$k\in \NN^p,\ l\in\NN^q,\  m\in \NN^r$, \\
$k_1\le k_2\le\dots\le k_p,\ l_1\le\dots\le l_q,\ m_1\le\dots\le m_r$, \\
$\rho=\rho_k^+\oplus\rho_l^-\oplus\rho'_m$, \\
$\alpha=0$,\\
$\gamma(H,X,Y)=c,\ c\in\RR$.
\end{enumerate}
\end{enumerate}
\end{theo}
\proof
We already know that for a given symmetric triple $(\fg,\theta,\ip)$ the Lie algebra 
$\fg/\fri(\fg)^\perp$ with the involution induced by $\theta$ is isomorphic to one of 
the Lie algebras with involution in Prop.~\ref{P71}. Moreover, if we consider 
$\fa:=\fri(\fg)^\perp/\fri(\fg)$ with the induced involution and the induced scalar 
product, then $\fa_-$ is positive definite if  $\dim \fl_-=2$ and the induced scalar 
product on $\fa_-$ has index $1$ if $\dim \fl_-=1$.

According to Theorem \ref{Pwieder} we have to determine the orbit space of the action of 
$G=\Aut(\fl,\thl)\times \Aut(\fa)$ on
    $$H:=\coprod_{\rho\in \Hom((\fl,\thl),\so(\fa))_{\rm ss}}
\cH^2_Q(\fl,\thl,\fa_\rho)_{0}\
$$
for these combinations of $(\fl,\thl)$ and $\fa$.

We begin with $\fl\in\{\fn(2), \,\fr_{3,-1}\}$. We fix an arbitrary 
pseudo-Euclidean space $\fa$ with involution $\tha$ such that $\fa_-$ is 
positive definite, i.e. $\fa=\fa_+^{\bar p, \bar q} \oplus \fa_-^{0, \bar s}$.  
Here we use the description of the orbit space as $$\coprod_{[\rho]\in 
\Hom((\fl,\thl),\so(\fa))_{\rm ss}/G}
\cH^2_Q(\fl,\thl,\fa_\rho)_{0}/G_\rho\ .$$
{}From Propositions \ref{PAut} and \ref{PAut2} we know that
\begin{eqnarray*}
\Hom((\fl,\thl),\so(\fa))_{\rm ss}/G&=&\Hom((\fl,\thl),\so
(\fa))_{\rm ss}/(\Aut(\fa)\times \ZZ_2)\\ &=&\Hom((\fl,\thl),\so
(\fa))_{\rm ss}/\Aut(\fa)
\end{eqnarray*}
and by Prop.~\ref{PHom} we can identify $\Hom((\fl,\thl),\so
(\fa))_{\rm ss}/\Aut(\fa)$ with
$$\{\ (\rho_0)^{p_0,q_0}_{0,s_0} \oplus \tilde\rho^+_\lambda \oplus \rho_\mu^+ 
\mid \lambda\in\Lambda_p,\,\mu\in\Lambda_q,\, p_0+p =\bar p,\ q_0+q =\bar q,\ 
p+q+s_0=\bar s \} \,.$$

Furthermore, we identify $\lambda\in(\fl_0^*)^p$ with $\lambda(X)\in\RR^p$ and 
therefore $\Lambda_p$ with 
$$\{\lambda=(\lambda^1,\dots,\lambda^p)\in\RR^p\mid 0< \lambda^1\le \dots
\le\lambda^p\}.$$

Now let $\rho=(\rho_0)^{p_0,q_0}_{0,s_0} \oplus \tilde\rho^+_\lambda \oplus 
\rho_\mu^+$ be fixed. Recall that we know $\cH^2_Q(\fl,\thl,\fa_\rho)_{0}$
(Prop.~\ref{Pn2}) and $\Aut(\fl,\thl)$ (Prop.~\ref{PAut}). Consider
$\alpha_1,\,\alpha_2\in Z_\fl$. Then $[\alpha_1,0]$ and $[\alpha_2,0]$ are on 
the same $G_\rho$-orbit if and only if
there exists a real number $a>0$ such that the conditions
$$ \alpha_1(Y,Z)=\pm a \alpha_2(Y,Z) $$ and 
$$\langle \alpha_1(X,Y), \alpha_1(X,Y)\rangle = a\langle \alpha_2(X,Y), 
\alpha_2(X,Y)\rangle \ \mbox{ if }\ \fl=\fn(2)$$
$$\langle \alpha_1(X,Y), \alpha_1(X,Z)\rangle = a\langle \alpha_2(X,Y), 
\alpha_2(X,Z)\rangle\ \mbox{ if }\ \fl=\fr_{3,-1}$$
are satisfied. Furthermore, $(C^3(\fl)\setminus 0)/G_\rho = (C^3(\fl)\setminus
0)/\Aut(\fl,\thl)=\ZZ_2$. This yields the assertion for $\fl\in\{\fn(2), 
\,\fr_{3,-1}\}$.

Now we consider $\fl=\fh(1)$. We define subsets $H^{p,q}_i\subset H$, $i=1,2$, 
$p,q\ge 0$ by
$$ H^{p,q}_i:=\{ [\alpha_i,0] \in \cH^2_Q(\fl,\thl,\fa_\rho)_{0}\mid
\rho= (\rho_0)^{0,0}_{0,i} \oplus \tilde\rho^+_\lambda \oplus \rho_\mu^+,\ 
\lambda\in (\fl_0^*\setminus 0)^p,\, \mu\in (\fl_0^*\setminus 0)^q \},$$
where $\alpha_1$ is given by 
$\alpha_1(X,Z)=A_0, \ \alpha_1(X,Y)=\alpha_1(Y,Z)=0$ and $\alpha_2$ is given by
$\alpha_2(X,Z)=A_1, \ \alpha_2(Y,Z)=A_2,\ \alpha_2(X,Y)=0$
for standard bases $A_0$ and $A_1,\,A_2$ of $\fa^\fl_\rho$, respectively. Note 
that $\alpha_1,\,\alpha_2\in Z_\fl$ for all such $\rho$. By our description of 
$\Aut(\fl,\thl)$ (Prop.~\ref{PAut2}) we know that for each orbit $\cal O$ of the 
$G$-action on $H$ the intersection ${\cal O}\cap H^{p,q}_i$ is non-empty for 
exactly one triple $(i,p,q)$. Moreover, ${\cal O}\cap H^{p,q}_i$ is a $G_i\times 
\bar \frak S_{p,q}$-orbit in $H^{p,q}_i$, where 
$G_i=(\Aut(\fl,\thl)\times\Aut(\fa^\fl))_{\alpha_i}$ is the stabilizer of 
$\alpha_i\in C^2(\fl,\fa^\fl)$. Using again our description of $\Aut(\fl,\thl)$
from Prop.~\ref{PAut2} it is not hard to compute $G_i$ and then $H^{p,q}_i/G_i$.

In the cases $\fl=\RR^1$ and $\fl=\RR^2$ we proceed as in the case of $\fl=\fh(1)$.

The result for $\fl\in\{\fsu(2),\fsl(2,\RR)\}$ is a straightforward consequence of
Lemma \ref{rep} and Lemma \ref{whitehead}. Note that in this case any automorphism of 
$\fl$ acts trivially
on $C^3(\fl)$.
\qed

Recall that a {\em pseudo-Hermitian symmetric space} is a pseudo-Riemannian symmetric 
space
$(M,g)$ equipped with an almost complex structure $J$ which is compatible with $g$ and 
the involutions $\theta_x$, $x\in M$. Then $J$ is automatically integrable and $(M,g,J)$
is a K\"ahler manifold (see \cite{CP}, Ch.I,\S 6).

\begin{co}\label{kahl}
If $(\fg,\theta,\ip)$ is a symmetric triple associated
with a simply connected indecomposable pseudo-Hermitian symmetric
space of complex signature $(1,q)$, $q\ge 0$, which is neither semi-simple
nor flat,
then $(\fg,\theta,\ip)$ is isomorphic to
$\dd$ for exactly one of the data in the following list (which contains only data giving 
rise to such triples):
\begin{enumerate}
\item $q=1:$ $\fl=\fl_-=\RR^2=\Span\{Y,Z\}$,
\begin{enumerate}
\item $\fa=\fa_+^{0,1}$, $\rho=(\rho_0)^{0,1}_{0,0}$, $\alpha(Y,Z)=A_0$, where $A_0$ is 
the
standard basis of $\fa_+^{0,1}$,
$\gamma=0$;
\item $\fa=\fa_+^{1,0}$, $\rho=(\rho_0)^{1,0}_{0,0}$, $\alpha(Y,Z)=A_0$, where $A_0$ is 
the
standard basis of $\fa_+^{1,0}$,
$\gamma=0$;
\end{enumerate}
\item $q=2:$ $\fl=\fh(1),\ \fl_+=\RR\cdot Z ,\ \fl_-=\Span\{X,Y \} $
\begin{enumerate}
\item[] $\fa=\fa_-^{0,2}$, $\rho=(\rho_0)^{0,0}_{0,2}$, $\alpha(X,Z)=A_1$, 
$\alpha(Y,Z)=A_2$,
$\alpha(X,Y)=0$,\\ where $A_1,A_2$ is the standard basis of $\fa_-^{0,2}$,\\
$\gamma=0$;
\end{enumerate}
\item $q=1+p,\,p\ge 0:$ $\fl=\fsu(2),\ \fl_+= \RR\cdot H,\ \fl_-=\Span\{X,\ Y\},$
\begin{enumerate}
\item[] $\fa=\fa_+^{0,2p-r}\oplus\fa_-^{0,2p}$, $\ 0\le r\le p$,
$\ \displaystyle
\rho=\bigoplus_{i=1}^r\rho_1^-\oplus\bigoplus_{i=1}^{p-r}\rho'_1$,\\
$\alpha=0$,
$\gamma(H,X,Y)=c,\ c\in\RR$.
\end{enumerate}
\item $q=1+p,\,p\ge 0:$ $\fl=\fsl(2,\RR), \ \fl_+=\RR\cdot (X-Y),\
\fl_-=\Span\{H,\ X+Y \}.$
\begin{enumerate}
\item[] $\fa=\fa_+^{2p-r,0}\oplus\fa_-^{0,2p}$, $\ 0\le r\le p$,
$\ \displaystyle
\rho=\bigoplus_{i=1}^r\rho_1^-\oplus\bigoplus_{i=1}^{p-r}\rho'_1$,\\
$\alpha=0$,
$\gamma(H,X,Y)=c,\ c\in\RR$.
\end{enumerate}
\end{enumerate}
\end{co}
\proof A symmetric triple $(\fg,\theta,\ip)$ associated with a simply connected non-flat 
indecomposable pseudo-Hermitian symmetric
space $M$ is indecomposable (\cite{CP}, Ch.I, Prop. 6.6). Moreover, it carries an 
automorphism
$F$ satisfying $F^2=\theta$, $F|_{\fg_+}=\Id$, and any such automorphism induces a 
pseudo-Hermitian structure $J$ on $M$ (\cite{CP}, Ch.I, Prop. 6.4).
If $(\fg,\theta,\ip)=\dd$ for $[\alpha,\gamma]\in \cH^2_Q(\fl,\thl,\fa)_\sharp$,
then the existence of such an automorphism $F$ is equivalent to the existence
of an automorphism $\bar F=(F_\fl,F_\fa^{-1})$ of the triple $(\fl,\thl,\fa)$
satisfying
\begin{equation}\label{fanta}
F_\fl^{\, 2}=\theta_\fl,\ F_\fl |_{\fl_+}=\Id,\quad F_\fa^{\, 2}=\theta_\fa,\  
F|_{\fa_+}=\Id
\end{equation}
such that $(\bar F^*\alpha,\bar F^*\gamma)=(\alpha,\gamma)$. The corollary now follows
by inspection of the classification list in Theorem \ref{Klasse} provided
one takes the following fact into account:

Let $(\rho,\fa)$ be a semi-simple orthogonal $(\fl,\thl)$-module
such that the triple $(\fl,\thl,\fa)$ admits an automorphism $\bar F$ satisfying 
(\ref{fanta}).
\begin{enumerate}
\item[(i)](\cite{CP}, Ch.V, Prop.3.3) If $(\fl,\thl)$ is as in cases 6. or 8. of
Theorem \ref{Klasse} and $\fa_-$ is positive definite, then $\rho$ is the direct sum of 
irreducible $(\fl,\thl)$-modules equivalent to $\rho_1^-$ or $\rho_1'$. All these 
modules
admit an automorphism $F_\fa$ with the required properties.
\item[(ii)] If $\fl$ is solvable and satisfies $(T_1)$, then $\rho$ is the trivial 
representation on $\fa$.
\end{enumerate}

It remains to prove (ii). Let $\lambda\in \fl_0^*\setminus 0$ and let
$E_\lambda\subset \fa_{\Bbb C}$ be the
corresponding weight space. We have to show that $E_\lambda=\{0\}$. The automorphism
$F_\fl$ induces a complex structure $j$ on the real vector space $\fl_0^*$. The elements 
$j^k(\lambda)$, $k=0,1,2,3$, are pairwise different. Therefore the sum of the weight 
spaces $E_{j^k(\lambda)}$, $k=0,1,2,3$, is direct. Take $v\in E_\lambda$.
Then $F_\fa^{\, k}(v)\in E_{j^k(\lambda)}$,
and $v^-:=v-F_\fa(v)+F_\fa^{\, 2}(v)-F_\fa^{\, 3}(v)$ satisfies $F_\fa(v^-)=-v^-$. 
However, the only
eigenvalues of $F_\fa$ on $\fa_{\Bbb C}$ are $1$, $i$ and $-i$. We conclude that 
$v^-=0$, hence $v=0$.
This finishes the proof of (ii).
\qed

Note that Assertion (ii) of the above proof has strong implications for the structure
of arbitrary pseudo-Hermitian symmetric spaces with solvable transvection group.

\vspace{0.8cm}
{\footnotesize
Ines Kath\\
Max-Planck-Institut f{\"u}r Mathematik in den Naturwissenschaften\\
Inselstr. 22-26, D-04103 Leipzig, Germany\\
email: ikath@mis.mpg.de\\[2ex]
Martin Olbrich\\
Mathematisches Institut der Georg-August-Universit{\"a}t G{\"o}ttingen\\
Bunsenstr. 3-5, D-37073 G{\"o}ttingen, Germany\\
email: olbrich@uni-math.gwdg.de
}

\begin{thebibliography}{MMMM}

\bibitem[ABCV]{ABCV} D.\,V.\,Alekseevsky, N.\,Bla\v zi\'c, V.\,Cort\'es, S.\,Vukmirovi\'c,
{\sl A class of Osserman spaces.} Preprint no. 3, SFB 611, University Bonn, 10 pages.

\bibitem[AC\,1]{AC1} D.\,V.\,Alekseevsky, V.\,Cort\'es,
{\sl Classification of indefinite hyper-K\"ahler symmetric spaces.}
Asian J.\,Math. {\bf 5} (2001), 663 -- 684.

\bibitem[AC\,2]{AC2} D.\,V.\,Alekseevsky, V.\,Cort\'es,
{\sl Classification of pseudo-Riemannian symmetric spaces of quaternionic K\"ahler type.} To appear in Amer.\,Math.\,Soc.\,Transl.

\bibitem[B]{B} M.\,Berger,
{\sl Les espaces sym\'etriques non compacts.}
Ann.\,Sci.\,Ecole\,Norm. Sup. {\bf 74} (1957), 85 -- 177.

\bibitem[BB\,1]{BB1}
L.\,Berard Bergery,
{\sl D\'ecomposition de Jordan-H\"older d'une repr\'esentation de dimension finie,
adapt\'ee \`a une forme r\'eflexive.} Handwritten notes.

\bibitem[BB\,2]{BB2}
L.\,Berard Bergery,
{\sl Structure des espaces sym\'etriques pseudo-riemanniens.}
Handwritten notes.

\bibitem[CP\,1]{CP1} M.\,Cahen, M.\,Parker, {\sl Sur des classes d'espaces
pseudo-riemanniens sym\'et\-riques.}
Bull.\,Soc.\,Math.\,Belg. {\bf 22} (1970), 339 -- 354.

\bibitem[CP\,2]{CP} M.\,Cahen, M.\,Parker, {\sl Pseudo-Riemannian symmetric spaces.}
Mem.\,Amer. Math.\,Soc. {\bf 24} (1980), no. 229.

\bibitem[CW]{CW} M.\,Cahen, N.\,Wallach, {\sl Lorentzian symmetric spaces.}
Bull.\,Amer.\,Math.\,Soc.  {\bf 76} (1970), 585--591.

\bibitem[H]{H} S.\,Helgason, {\sl Differential Geometry, Lie Groups, and
Symmetric Spaces.} Graduate Studies in Mathematics 34,
AMS, Providence, RI, 2001.

\bibitem[KO\,1]{KO1} I.\,Kath, M.\,Olbrich,
{\sl Metric Lie algebras with maximal isotropic centre.}
Math. Z. {\bf 246} (2004), 23 -- 53.

\bibitem[KO\,2]{KO2} I.\,Kath, M.\,Olbrich,
{\sl Metric Lie algebras and quadratic extensions.}
arXiv: math.DG/0312243.

\bibitem[KN]{KN} S.\,Kobayashi, K.\,Nomizu, {\sl Foundations of differential
geometry II.}  New York, 1996.

\bibitem[N]{N}Th.\,Neukirchner, {\sl Solvable Pseudo-Riemannian
    Symmetric Spaces.} arXiv: math.DG/0301326.
\end{thebibliography}
\end{document}